\documentclass[12pt,twoside]{article}
\usepackage[english]{babel}
\usepackage[latin1]{inputenc}
\usepackage{amsmath}
\usepackage{amssymb,amsfonts}
\usepackage{graphicx}
\usepackage{times,amssymb,amscd}
\newtheorem{thm}{Theorem}[section]

\newtheorem{lem}[thm]{Lemma}

\newtheorem{defn}[thm]{Definition}
\newtheorem{rem}[thm]{Remark}
\numberwithin{equation}{section}

\newcommand{\bA}{\mathbf{A}}
\newcommand{\bC}{\mathbf{C}}
\newcommand{\bE}{\mathbf{E}}
\newcommand{\bG}{\mathbf{G}}
\newcommand{\bH}{\mathbf{H}}

\newcommand{\bL}{\mathbf{L}}

\newcommand{\bR}{\mathbf{R}}
\newcommand{\bS}{\mathbf{S}}
\newcommand{\bV}{\mathbf{V}}

\newcommand{\be}{\mathbf{e}}

\newcommand{\br}{\mathbf{r}}
\newcommand{\bx}{\mathbf{x}}

\newcommand{\bT}{\mathbf{T}}

\newcommand{\bp}{\mathbf{p}}
\newcommand{\bt}{\mathbf{t}}

\newcommand{\BV}{\boldsymbol{V}}

\newcommand{\Be}{\boldsymbol{e}}
\newcommand{\Bu}{\boldsymbol{u}}
\newcommand{\Bv}{\boldsymbol{v}}

\newcommand{\cP}{\mathcal{P}}
\newcommand{\cS}{\mathcal{S}}

\newcommand{\cH}{\mathcal{H}}

\newcommand{\cA}{\mathcal{A}}

\newcommand{\EUC}{\mathbf E^3}
\newcommand{\SPH}{\bS^3}
\newcommand{\HYP}{\bH^3}
\newcommand{\SXR}{\bS^2\!\times\!\mathbb{R}}
\newcommand{\HXR}{\bH^2\!\times\!\mathbb{R}}
\newcommand{\SLR}{\widetilde{\bS\bL_2\bR}}
\newcommand{\NIL}{\mathbf{Nil}}
\newcommand{\SOL}{\mathbf{Sol}}

\begin{document}
\pagestyle{myheadings}
\markboth{\centerline{G\'eza Csima and Jen\H o Szirmai}}
{Translation-like Apollonius and triangular surfaces  $\dots$}
\title
{Translation-like Apollonius and triangular surfaces in non-constant curvature Thurston geometries
}

\author{G\'eza Csima and Jen\H o Szirmai\\
\normalsize Department of Algebra and Geometry, Institute of Mathematics,\\
\normalsize Budapest University of Technology and Economics, \\
\normalsize M\"uegyetem rkp. 3., H-1111 Budapest, Hungary \\
\normalsize csimageza@gmail.com, szirmai@math.bme.hu
\date{\normalsize{\today}}}

\maketitle
\begin{abstract}
In the present paper we deal with non-constant curvature Thurston geometries \cite{M97}, \cite{S}, \cite{Sz22-3},\cite{W06}.
We define and determine the generalized trans\-lation-like Apollonius surfaces and thus also bisector surfaces as a special case.
Moreover, we give a possible definition of the ``surface of a translation-like triangle" in each investigated geometry.
In our work we will use the projective model of Thurston geometries described by E. Moln\'ar in \cite{M97}.
\end{abstract}

\section{Introduction, preliminary results} \label{section1}
In classical geometries with constant curvature, Dirichlet-Voronoi cells (brifly D-V cells) and their tilings of geometric spaces play a fundamental role. 
Another important aspect is related to material structure issues, since they also play a fundamental role in crystallographyy. 

D-V cells can be derived using bisector surfaces, which are a special case of Apollonius surfaces. Apollonius surfaces and related theorems 
(such as Ceva and Menealaus theorems \cite{Sz21}, \cite{Sz22}) play an important role in the internal structure of geometries. 

In the Thurston spaces, one can introduce in a natural way (see \cite{M97}) 
translation mappings any point to any other point. Consider a unit tangent 
vector at the origin. Translations carry this vector to a tangent vector any 
other point.
If a curve $t\rightarrow (x(t),y(t),z(t))$ has just the translated vector as its tangent vector at
each point, then the curve is called a translation curve. This assumption leads to
a system of first order differential equations. Thus translation curves are simpler
than geodesics and differ from them in $\NIL$, $\SLR$ and $\SOL$ geometries. 
In $\EUC$, $\SPH$, $\HYP$, $\SXR$ and $\HXR$ geometries, the translation and geodesic curves coincide
with each other. But in the $\NIL$, $\SLR$ and $\SOL$ geometries, 
translation curves are in many ways more natural than geodesics. Therefore, we
distinguish two different distance functions: $d^g$ is the usual geodesic distance function,
and $d^t$ is the translation distance function. So we obtain
two types of the Apollonius or bisector surfaces (and two types of the corresponding D-V cells) from the two different distance functions, 
but {\it{in the present paper we consider only the translation case}}.

The classical definition of the Apollonius circle in the Euclidean plane $\mathbf{E}^2$ is  
the set of all points of $\mathbf{E}^2$ whose distances from two fixed points are in a constant ratio $\lambda\in\mathbf{R}^+$. In \cite{Sz21} we extended this definition 
in a natural way to the further Thurston geometries 
$$
\EUC,\SPH,\HYP,\SXR,\HXR,\NIL,\SLR,\SOL.
$$
\begin{defn}
The Apollonius surface in the Thurston geometry $X \in \{\EUC,\SPH,\HYP, $ $\SXR,\HXR,\NIL,\SLR,\SOL\}$ is  
the set of all points of $X$ whose translation distances from two fixed points are in a constant ratio $\sigma\in\mathbf{R}^+$.
\end{defn}
A special case of Apollonius surfaces is the translation-like bisector (or equidistant) surface ($\lambda=1)$
of two arbitrary points of $X$. 
In the Thurston geometries $\SXR$, $\HXR$,
$\NIL$, $\SOL$, $\SLR$ there are few results on this topic. Let $X$ be one of the above five geometries and let $\Gamma$ be a discrete group of isometries of $X$ and $d^t$ is the translation distance function. 

We define the Dirichlet-Voronoi cell with kernel point $K$ of a given discrete isometry group $\Gamma$:
\begin{defn}
We say that the point set
\[
\mathcal{D}(K)=\left\{ Y \in X : ~ d^t(K,Y)\le d^t(K^{\mathbf{g}},Y) ~ \text{for all} ~ \mathbf{g}\in \Gamma \right\}\subset X
\]
is the {\it Dirichlet-Voronoi cell} of $\Gamma$ around its {\it kernel point} $K$ where $d^t$ is the translation distance function of $X$.
\end{defn}
These Apollonius and bisector surfaces have an important role in structure of Dirichlet - Voronoi (briefly, D-V) cells and they are relevant in the study of tilings, ball packing and ball covering. In $3$-dimensional spaces of constant curvature, the
D-V cells have been widely investigated, but in the other Thurston geometries $\SXR$, $\HXR$,
$\NIL$, $\SOL$, $\SLR$ there are few results on this topic. 

In \cite{PSSz10}, \cite{PSSz11-1}, \cite{PSSz11-2} we studied the geodesic-like bisector surfaces and 
the Apollonius surfaces are investigated in $\SXR$, $\HXR$ and $\NIL$ geometries in \cite{Sz21, Sz22}. Moreover, we studied in \cite{Sz19} and \cite{VSz19} the translation-like equidistant surfaces in $\SOL$ and $\NIL$ geometries.

The next important question, which is also closely related to the 
Apollonius surfaces, is the determination of the surface of the given geodesic or translational triangle. Defining this is an essential condition for stating elementary geometric concepts and theorems related to triangles.

However, defining the surface of a geodesic triangle in $X$ space is not straightforward (see \cite{Sz21,Sz22}). The usual geodesic triangle surface definition is not possible because the geodesic or translation curves starting from different vertices and ending at points of the corresponding opposite edges define different 
surfaces in general, i.e. {\it geodesics or translation curves starting from different vertices and ending at points on the corresponding opposite side usually do not intersect.}

In the works, \cite{Sz21}, \cite{Sz22} we proposed a possible definition of a surface in $\SXR$, $\HXR$ and $\NIL$ geometries, which was done with the help of Apollonius surfaces and returned traditional triangular surfaces in geometries with constant curvature.
 
{\it In this work, we extend these questions to the translation curves of non-constant curvature Thurston geometries.} We study the translation-like Apollonius surfaces of two points in $\NIL$, $\SOL$ and $\SLR$ geometries, determine their equations and visualize them. Moreover, we also provide a new possible definition of the surfaces of translation-like triangles.
\section{The structures of the investigated geometries}
In geometries with constant curvature $\EUC$, $\HYP$, $\SPH$ 
the translation curves and the geodesic curves are the same, and in these the Apollonius surfaces and the surfaces of the geodesic triangles are well defined. 


In the following, we summarize the most important concepts of the further investigated geometries.
\subsection{On $\SXR$ and $\HXR$ geometries}
In \cite{M97} E.~Moln\'ar has shown that the homogeneous 3-spaces have a unified interpretation in the projective 3-sphere $\mathcal{P} \mathcal{S}^3(\BV^4,\BV_4,\mathbf{R}).$ In this work, we will use this projective model of $\SXR$ and $\HXR.$
We will use the Cartesian homogeneous coordinate simplex $E_0(\be_0)$,$E_1^{\infty}(\be_1)$,$E_2^{\infty}(\be_2)$,
$E_3^{\infty}(\be_3),$ $(\{\be_i\}\subset \bV^4$ with the unit point $E(\be = \be_0 + \be_1 + \be_2 + \be_3 ))$ 
which is distinguished by an origin $E_0$ and by the ideal points of coordinate axes, respectively. 
Moreover, $\mathbf{y}=c\bx$ with $0<c\in \mathbb{R}$ (or $c\in\mathbf{R}\setminus\{0\})$
defines a point $(\bx)=(\mathbf{y})$ of the projective 3-sphere $\mathcal{P} \mathcal{S}^3$ (or that of the projective space $\cP^3$ where opposite rays
$(\bx)$ and $(-\bx)$ are identified). 
The dual system $\{(\Be^i)\}, \ (\{\Be^i\}\subset \BV_4)$, with $\be_i\Be^j=\delta_i^j$ (the Kronecker symbol), describes the simplex planes, especially the plane at infinity 
$(\Be^0)=E_1^{\infty}E_2^{\infty}E_3^{\infty}$, and generally, $\Bv=\Bu\frac{1}{c}$ defines a plane $(\Bu)=(\Bv)$ of $\cP \cS^3$
(or that of $\cP^3$). Thus $0=\bx\Bu=\mathbf{y}\Bv$ defines the incidence of point $(\bx)=(\mathbf{y})$ and plane
$(\Bu)=(\Bv)$, as $(\bx) \text{I} (\Bu)$ also denotes it. Thus $\SXR$ can be visualized in the affine 3-space $\bA^3$
(so in $\bE^3$) and the points of $\HXR$ form an open cone solid, described by the following set:
\begin{equation}
\HXR=\left\{X(\bx=x^i\be_i)\in\mathcal{P}^3:-(x^1)^2+(x^2)^2+(x^3)^2<0<x^0,x^1\right\}
\label{2.1}
\end{equation}

%
In both geometries, we introduce a new coordinate system in order to write the translation curves more simply. We introduce the polar parametrization of $\SXR$ and the cylindrical parametrization of $\HXR$ in $\BV^4$:
\begin{equation}
\begin{gathered}
\SXR:\  x_0=1,\ x_1 = e^t\cos\theta,\ x_2 =e^t \sin\theta\cos\phi,\ x_3 =e^t\sin\theta\sin\phi \\
-\pi<\phi\leq\pi,\ 0\leq\theta\leq\pi,\ t\in\mathbf{R}
\label{2.3}
\end{gathered}
\end{equation}
\begin{equation}
\begin{gathered}
\HXR:\  x_0=1,\ x_1 = e^t\cosh r,\ x_2 =e^t \sinh r\cos\phi,\ x_3 =e^t\sinh r\cos\phi \\
-\pi<\phi\leq\pi,\ 0\leq r,\ t\in\mathbf{R}
\label{2.4}
\end{gathered}
\end{equation}
where $(\theta,\phi)$ and $(r,\phi)$ are the usual polar coordinates of $\bS^2$ and $\bH^2,$ furthermore $t$ is the real component, the so-called fibre coordinate in the direct product $\SXR$ and $\HXR.$ With $x = \frac{x_1}{x_0},$ $y = \frac{x_2}{x_0},$ $z = \frac{x_3}{x_0},$ setting $t$ to be $0$ describes the unit sphere in \ref{2.3}, and the $x>0$ part of the two-sheeted hyperboloid $x^2-y^2-z^2=1$ in \ref{2.4}. This last surface can be called the unit hyperboloid of $\HXR.$ In both geometries $t=\infty$ would be the ideal plane $(\be^0)$ at infinity, $t = -\infty$ would be the origin $(\be_0)$ in limit in $\EUC$ model. Central similarity with factor $e^a$ means the translation by $\mathbf{R}$-component $a,$ commuting with any isometry of $\bS^2$ and $\bH^2.$


\subsubsection{Translation curves}

We can assume that the starting point of a translation curve in both geometries is the $(1,1,0,0)$ point. Hereafter, let the functions $S(t)$ and $C(t)$ be $\sin(t)$ and $\cos(t)$ in $\SXR,$ $\sinh(t)$ and $\cosh(t)$ in $\HXR.$ Then the translation curve by \cite{MoSzi10} can be given:
\begin{equation}
\begin{gathered}
x(\tau)=e^{\tau\sin (v)}C(\tau\cos (v)),\\
y(\tau)=e^{\tau\sin (v)}S(\tau\cos (v))\cos(u),\\
z(\tau)=e^{\tau\sin (v)}S(\tau\cos (v))\sin(u),\\
-\pi<u\leq\pi,\ -\frac{\pi}{2}\leq v\leq\frac{\pi}{2}.
\label{2.6}
\end{gathered}
\end{equation}

In the parametric equation of the translation curve above $\tau$ denotes the arc-length parameter; $v$ denotes the angle, formed by the starting direction vector of the curve and the tangent plane at the origin $E_1=(1,1,0,0)$ of the unit sphere ($x^2+y^2+z^2=1$) for $\SXR$ and the tangent plane at $E_1$ of unit hyperboloid ($x^2-y^2-z^2=1$) for $\HXR$; and $u$ denotes the angle, formed by the $y$ axis and the projection of this starting direction onto the tangent plane, described above. 
\begin{rem}
\begin{enumerate}
	\item It is easy to see, that the translation curve lies in a plane with equation:
\begin{equation}
\sin(u)y-\cos(u) z=0  
\label{2.5}
\end{equation}
\item The tangent vector of (\ref{2.6}) at the origin is provided by replacing 0 for $\tau$ in the derivative by $\tau.$
\begin{equation}
\bt=(\sin(v),\cos(v)\cos(u),\cos(v)\sin(u)) 
\label{2.6a}
\end{equation}
\end{enumerate}
\end{rem}
\subsection{On $\NIL$ geometry}
$\NIL$ geometry can be derived from the famous real matrix group $\mathbf{L(\mathbf{R})}$ discovered by Werner Heisenberg. The left (row-column) 
multiplication of Heisenberg matrices
     \begin{equation}
     \begin{gathered}
     \begin{pmatrix}
         1&x&z \\
         0&1&y \\
         0&0&1 \\
       \end{pmatrix}
       \begin{pmatrix}
         1&a&c \\
         0&1&b \\
         0&0&1 \\
       \end{pmatrix}
       =\begin{pmatrix}
         1&a+x&c+xb+z \\
         0&1&b+y \\
         0&0&1 \\
       \end{pmatrix}
      \end{gathered} \label{2.7}
     \end{equation}
defines "translations" $\mathbf{L}(\mathbf{R})= \{(x,y,z): x,~y,~z\in \mathbf{R} \}$ 
on the points of $\NIL= \{(a,b,c):a,~b,~c \in \mathbb{R}\}$. 
These translations are not commutative in general. The matrices $\mathbf{K}(z) \vartriangleleft \mathbf{L}$ of the form
     \begin{equation}
     \begin{gathered}
       \mathbf{K}(z) \ni
       \begin{pmatrix}
         1&0&z \\
         0&1&0 \\
         0&0&1 \\
       \end{pmatrix}
       \mapsto (0,0,z)  
      \end{gathered}\label{2.8}
     \end{equation} 
constitute the one parametric centre, i.e. each of its elements commutes with all elements of $\mathbf{L}$. 
The elements of $\mathbf{K}$ are called {\it fibre translations}. $\NIL$ geometry of the Heisenberg group can be projectively 
(affinely) interpreted by "right translations" 
on points as the matrix formula 
     \begin{equation}
     \begin{gathered}
       (1;a,b,c) \to (1;a,b,c)
       \begin{pmatrix}
         1&x&y&z \\
         0&1&0&0 \\
         0&0&1&x \\
         0&0&0&1 \\
       \end{pmatrix}
       =(1;x+a,y+b,z+bx+c) 
      \end{gathered} \label{2.9}
     \end{equation} 
shows, according to (\ref{2.7}). Here we consider $\mathbf{L}$ as projective collineation group with right actions in homogeneous coordinates.
We will use the usual projective model of $\NIL$ (see \cite{M97} and \cite{MSz06}).

The translation group $\mathbf{L}$ defined by formula (\ref{2.9}) can be extended to a larger group $\mathbf{G}$ of collineations,
preserving the fibres, that will be equivalent to the (orientation preserving) isometry group of $\NIL$. 

In \cite{MSz06} we has shown that 
a rotation through angle $\omega$
about the $z$-axis at the origin, as isometry of $\NIL$, keeping invariant the Riemann
metric everywhere, will be a quadratic mapping in $x,y$ to $z$-image $\overline{z}$ as follows:
     \begin{equation}
     \begin{gathered}
       \mathcal{M}=\br(O,\omega):(1;x,y,z) \to (1;\overline{x},\overline{y},\overline{z}); \\ 
       \overline{x}=x\cos{\omega}-y\sin{\omega}, \ \ \overline{y}=x\sin{\omega}+y\cos{\omega}, \\
       \overline{z}=z-\frac{1}{2}xy+\frac{1}{4}(x^2-y^2)\sin{2\omega}+\frac{1}{2}xy\cos{2\omega}.
      \end{gathered} \label{2.10}
     \end{equation}
This rotation formula $\mathcal{M}$, however, is conjugate by the quadratic mapping $\alpha$ to the linear rotation $\Omega$ as follows
     \begin{equation}
     \begin{gathered}
       \alpha^{-1}: \ \ (1;x,y,z) \stackrel{\alpha^{-1}}{\longrightarrow} (1; x',y',z')=(1;x,y,z-\frac{1}{2}xy) \ \ \text{to} \\
       \Omega: \ \ (1;x',y',z') \stackrel{\Omega}{\longrightarrow} (1;x",y",z")=(1;x',y',z')
       \begin{pmatrix}
         1&0&0&0 \\
         0&\cos{\omega}&\sin{\omega}&0 \\
         0&-\sin{\omega}&\cos{\omega}&0 \\
         0&0&0&1 \\
       \end{pmatrix}, \\
       \text{with} \ \ \alpha: (1;x",y",z") \stackrel{\alpha}{\longrightarrow}  (1; \overline{x}, \overline{y},\overline{z})=(1; x",y",z"+\frac{1}{2}x"y").
      \end{gathered} \label{2.11}
     \end{equation}
This quadratic conjugacy modifies the $\NIL$ translations in (\ref{2.9}), as well. 
\subsubsection{Translation curves}
We consider a $\NIL$ curve $(1,x(t), y(t), z(t) )$ with a given starting tangent vector at the origin $O=E_0=(1,0,0,0)$
\begin{equation}
   \begin{gathered}
      u=\dot{x}(0),\ v=\dot{y}(0), \ w=\dot{z}(0).
       \end{gathered} \label{2.12}
     \end{equation}
For a translation curve let its tangent vector at the point $(1,x(t), y(t), z(t) )$ be defined by the matrix (\ref{2.9}) 
with the following equation:
\begin{equation}
     \begin{gathered}
     (0,u,v,w)
     \begin{pmatrix}
         1&x(t)&y(t)&z(t) \\
         0&1&0&0 \\
         0&0&1&x(t) \\
         0&0&0&1 \\
       \end{pmatrix}
       =(0,\dot{x}(t),\dot{y}(t),\dot{z}(t)).
       \end{gathered} \label{2.13}
     \end{equation}
Thus, the {\it translation curves} in $\NIL$ geometry (see  \cite{MSz06}, \cite{MoSzi10}, \cite{MSzV}) are defined by the above first order differential equation system 
$\dot{x}(t)=u, \ \dot{y}(t)=v,  \ \dot{z}(t)=v \cdot x(t)+w,$ whose solution is the following: 
\begin{equation}
   \begin{gathered}
       x(t)=u t, \ y(t)=v t,  \ z(t)=\frac{1}{2}uvt^2+wt.
       \end{gathered} \label{2.14}
\end{equation}
We assume that the starting point of a translation curve is the origin, because we can transform a curve into an 
arbitrary starting point by translation (\ref{2.9}), moreover, unit initial velocity translation can be assumed by "geographic" parameters $\phi$ and $\theta$:
\begin{equation}
\begin{gathered}
        x(0)=y(0)=z(0)=0; \\ \ u=\dot{x}(0)=\cos{\theta} \cos{\phi}, \ \ v=\dot{y}(0)=\cos{\theta} \sin{\phi}, \ \ w=\dot{z}(0)=\sin{\theta}; \\ 
        - \pi \leq \phi \leq \pi, \ -\frac{\pi}{2} \leq \theta \leq \frac{\pi}{2}. \label{2.15}
\end{gathered}
\end{equation}
\subsection{On Sol geometry}
\label{sec:1}
In this Section we summarize the significant notions and notations of real $\SOL$ geometry (see \cite{M97}, \cite{S}).

$\SOL$ is defined as a 3-dimensional real Lie group with multiplication
\begin{equation}
     \begin{gathered}
(a,b,c)(x,y,z)=(x + a e^{-z},y + b e^z ,z + c).
     \end{gathered} \label{2.16}
     \end{equation}
We note that the conjugacy by $(x,y,z)$ leaves invariant the plane $(a,b,c)$ with fixed $c$:
\begin{equation}
     \begin{gathered}
(x,y,z)^{-1}(a,b,c)(x,y,z)=(x(1-e^{-c})+a e^{-z},y(1-e^c)+b e^z ,c).
     \end{gathered} \label{2.17}
     \end{equation}
Moreover, for $c=0$, the action of $(x,y,z)$ is only by its $z$-component, where $(x,y,z)^{-1}=(-x e^{z}, -y e^{-z} ,-z)$. Thus the $(a,b,0)$ plane is distinguished as a {\it base plane} in
$\SOL$, or by other words, $(x,y,0)$ is normal subgroup of $\SOL$.
$\SOL$ multiplication can also be affinely (projectively) interpreted by ``right translations"
on its points as the following matrix formula shows, according to (\ref{2.16}):
     \begin{equation}
     \begin{gathered}
     (1,a,b,c) \to (1,a,b,c)
     \begin{pmatrix}
         1&x&y&z \\
         0&e^{-z}&0&0 \\
         0&0&e^z&0 \\
         0&0&0&1 \\
       \end{pmatrix}
       =(1,x + a e^{-z},y + b e^z ,z + c)
       \end{gathered} \label{2.18}
     \end{equation}
by row-column multiplication.
This defines ``translations" $\mathbf{L}(\mathbf{R})= \{(x,y,z): x,~y,~z\in \mathbf{R} \}$
on the points of space $\SOL= \{(a,b,c):a,~b,~c \in \mathbf{R}\}$.
These translations are not commutative, in general.
Here we can consider $\mathbf{L}$ as projective collineation group with right actions in homogeneous
coordinates as usual in classical affine-projective geometry.

We will use the usual projective model of $\SOL$ (see \cite{M97} and \cite{MSz06}).

It will be important for us that the full isometry group Isom$(\SOL)$ has eight components, since the stabilizer of the origin
is isomorphic to the dihedral group $\mathbf{D_4}$, generated by two involutive (involutory) transformations:
\begin{equation}
   \begin{gathered}
      (1)  \ \ y \leftrightarrow -y; \ \ (2)  \ x \leftrightarrow y; \ \ z \leftrightarrow -z; \ \ \text{i.e. first by $3\times 3$ matrices}:\\
     (1) \ \begin{pmatrix}
               1&0&0 \\
               0&-1&0 \\
               0&0&1 \\
     \end{pmatrix}; \ \ \
     (2) \ \begin{pmatrix}
               0&1&0 \\
               1&0&0 \\
               0&0&-1 \\
     \end{pmatrix}; \\
     \end{gathered} \label{2.19}
     \end{equation}
     with its product, generating a cyclic group $\mathbf{C_4}$ of order 4
     \begin{equation}
     \begin{gathered}
     \begin{pmatrix}
                    0&1&0 \\
                    -1&0&0 \\
                    0&0&-1 \\
     \end{pmatrix};\ \
     \begin{pmatrix}
               -1&0&0 \\
               0&-1&0 \\
               0&0&1 \\
     \end{pmatrix}; \ \
     \begin{pmatrix}
               0&-1&0 \\
               1&0&0 \\
               0&0&-1 \\
     \end{pmatrix};\ \
     \mathbf{Id}=\begin{pmatrix}
               1&0&0 \\
               0&1&0 \\
               0&0&1 \\
     \end{pmatrix}.
     \end{gathered} \notag
     \end{equation}
     Or we write by collineations fixing the origin $O=(1,0,0,0)$:
\begin{equation}
(1) \ \begin{pmatrix}
         1&0&0&0 \\
         0&1&0&0 \\
         0&0&-1&0 \\
         0&0&0&1 \\
       \end{pmatrix}, \ \
(2) \ \begin{pmatrix}
         1&0&0&0 \\
         0&0&1&0 \\
         0&1&0&0 \\
         0&0&0&-1 \\
       \end{pmatrix} \ \ \text{of form (2.19)}. \label{2.20}
\end{equation}
A general isometry of $\SOL$ to the origin $O$ is defined 
by a product $\gamma_O \tau_X$, first $\gamma_O$ of form (2.20) then 
$\tau_X$ of (\ref{2.18}). To
a general point $A=(1,a,b,c)$, this will be a product $\tau_A^{-1} \gamma_O \tau_X$, mapping $A$ into $X=(1,x,y,z)$.

We remark only that the role of $x$ and $y$ can be exchanged throughout the paper, but this leads to the mirror interpretation of $\SOL$.
As formula \ref{2.16} fixes the metric of $\SOL$, the change above is not an isometry of a fixed $\SOL$ interpretation. Other conventions are also accepted
and used in the literature.

{\it $\SOL$ is an affine metric space (affine-projective one in the sense of the unified formulation of \cite{M97}). Therefore, its linear, affine, unimodular,
etc. transformations are defined as those of the embedding affine space.}
\subsubsection{Translation curves}
We consider a $\SOL$ curve $(1,x(t), y(t), z(t) )$ with a given starting tangent vector at the origin $O=(1,0,0,0)$
\begin{equation}
   \begin{gathered}
      u=\dot{x}(0),\ v=\dot{y}(0), \ w=\dot{z}(0).
       \end{gathered} \label{2.21}
     \end{equation}
For a translation curve let its tangent vector at the point $(1,x(t), y(t), z(t) )$ be defined by the matrix (\ref{2.18})
with the following equation:
\begin{equation}
     \begin{gathered}
     (0,u,v,w)
     \begin{pmatrix}
         1&x(t)&y(t)&z(t) \\
         0&e^{-z(t)}&0&0 \\
         0&0&e^{z(t)}& 0 \\
         0&0&0&1 \\
       \end{pmatrix}
       =(0,\dot{x}(t),\dot{y}(t),\dot{z}(t)).
       \end{gathered} \label{2.22}
     \end{equation}
Thus, {\it translation curves} in $\SOL$ geometry (see \cite{MoSzi10} and \cite{MSz}) are defined by the first order differential equation system
$\dot{x}(t)=u e^{-z(t)}, \ \dot{y}(t)=v e^{z(t)},  \ \dot{z}(t)=w,$ whose solution is the following:
\begin{equation}
   \begin{gathered}
     x(t)=-\frac{u}{w} (e^{-wt}-1), \ y(t)=\frac{v}{w} (e^{wt}-1),  \ z(t)=wt, \ \mathrm{if} \ w \ne 0 \ \mathrm{and} \\
     x(t)=u t, \ y(t)=v t,  \ z(t)=z(0)=0 \ \ \mathrm{if} \ w =0.
       \end{gathered} \label{2.23}
\end{equation}
We assume that the starting point of a translation curve is the origin, because we can transform a curve into an
arbitrary starting point by translation (\ref{2.18}), moreover, unit velocity translation can be assumed :
\begin{equation}
\begin{gathered}
        x(0)=y(0)=z(0)=0; \\ \ u=\dot{x}(0)=\cos{\theta} \cos{\phi}, \ \ v=\dot{y}(0)=\cos{\theta} \sin{\phi}, \ \ w=\dot{z}(0)=\sin{\theta}; \\
        - \pi < \phi \leq \pi, \ -\frac{\pi}{2} \leq \theta \leq \frac{\pi}{2}. \label{2.24}
\end{gathered}
\end{equation}
Thus we obtain the parametric equation of the {\it translation curve segment} $t(\phi,\theta,t)$ with starting point at the origin in direction
\begin{equation}
\bt(\phi, \theta)=(\cos{\theta} \cos{\phi}, \cos{\theta} \sin{\phi}, \sin{\theta}) \label{2.25}
\end{equation}
where $t \in [0,r] ~ r \in \bR^+$. If $\theta \ne 0$ then the system of equation is:
\begin{equation}
\begin{gathered}
        \left\{ \begin{array}{ll}
        x(\phi,\theta,t)=-\cot{\theta} \cos{\phi} (e^{-t \sin{\theta}}-1), \\
        y(\phi,\theta,t)=\cot{\theta} \sin{\phi} (e^{t \sin{\theta}}-1), \\
        z(\phi,\theta,t)=t \sin{\theta}.
        \end{array} \right. \\
        \text{If $\theta=0$ then}: ~  x(t)=t\cos{\phi} , \ y(t)=t \sin{\phi},  \ z(t)=0.
        \label{2.26}
        \end{gathered}
\end{equation}
\subsection{On $\SLR$ geometry}
The real $ 2\times 2$ matrices $\begin{pmatrix}
         d&b \\
         c&a \\
         \end{pmatrix}$ with unit determinant $ad-bc=1$
constitute a Lie transformation group by the usual product operation, taken to act on row matrices as on point coordinates on the right as follows
\begin{equation}
\begin{gathered}
(z^0,z^1)\begin{pmatrix}
         d&b \\
         c&a \\
         \end{pmatrix}=(z^0d+z^1c, z^0 b+z^1a)=(w^0,w^1)\\
\mathrm{with} \ w=\frac{w^1}{w^0}=\frac{b+\frac{z^1}{z^0}a}{d+\frac{z^1}{z^0}c}=\frac{b+za}{d+zc}, \label{2.27}
\end{gathered}
\end{equation}
as action on the complex projective line $\bC^{\infty}$ (see \cite{M97}, \cite{MSz}).
This group is a $3$-dimensional manifold, because of its $3$ independent real coordinates and with its usual neighbourhood topology (\cite{S}, \cite{T}).
In order to model the above structure in the projective sphere $\cP \cS^3$ and in the projective space $\cP^3$ (see \cite{M97}),
we introduce the new projective coordinates $(x^0,x^1,x^2,x^3)$ where
$a:=x^0+x^3, \ b:=x^1+x^2, \ c:=-x^1+x^2, \ d:=x^0-x^3$
with the positive, then the non-zero multiplicative equivalence as projective freedom in $\cP \cS^3$ and in $\cP^3$, respectively.
Then it follows that $0>bc-ad=-x^0x^0-x^1x^1+x^2x^2+x^3x^3$
describes the interior of the above one-sheeted hyperboloid solid $\cH$ in the usual Euclidean coordinate simplex with the origin
$E_0(1;0;0;0)$ and the ideal points of the axes $E_1^\infty(0,1,0,0)$, $E_2^\infty(0,0,1,0)$, $E_3^\infty(0,0,0,1)$.
We consider the collineation group ${\bf G}_*$ that acts on the projective sphere $\cS\cP^3$  and preserves a polarity i.e. a scalar product of signature
$(- - + +)$, this group leaves the one sheeted hyperboloid solid $\cH$ invariant.
We have to choose an appropriate subgroup $\mathbf{G}$ of $\mathbf{G}_*$ as isometry group, then the universal covering group and space
$\widetilde{\cH}$ of $\cH$ will be the hyperboloid model of $\SLR$ \cite{M97}.

The elements of the isometry group of
$\mathbf{SL_2R}$ (and so by the above extension the isometries of $\SLR$) can be described in \cite{M97} and \cite{MSz}).
Moreover, we have the projective proportionality, of course.
We define the {\it translation group} $\bG_T$, as a subgroup of the isometry group of $\mathbf{SL_2R}$,
the isometries acting transitively on the points of ${\cH}$ and by the above extension on the points of $\SLR$ and $\widetilde{\cH}$.
$\bG_T$ maps the origin $E_0(1,0,0,0)$ onto $X(x^0,x^1,x^2,x^3)$. These isometries and their inverses (up to a positive determinant factor)
are given by the following matrices:
\begin{equation}
\begin{gathered} \bT:~(t_i^j)=
\begin{pmatrix}
x^0&x^1&x^2&x^3 \\
-x^1&x^0&x^3&-x^2 \\
x^2&x^3&x^0&x^1 \\
x^3&-x^2&-x^1&x^0
\end{pmatrix}.
\end{gathered} \label{2.28}
\end{equation}
The rotation about the fibre line through the origin $E_0(1;0;0;0)$ by angle $\omega$ $(-\pi<\omega\le \pi)$ can be expressed by the following matrix
(see \cite{M97})
\begin{equation}
\begin{gathered} \bR_{E_O}(\omega):~(r_i^j(E_0,\omega))=
\begin{pmatrix}
1&0&0&0 \\
0&1&0&0 \\
0&0&\cos{\omega}&\sin{\omega} \\
0&0&-\sin{\omega}&\cos{\omega}
\end{pmatrix},
\end{gathered} \label{2.29}
\end{equation}
and the rotation $\bR_X(\omega)$ about the fibre line through $X(x^0;x^1;x^2;x^3)$ by angle $\omega$ can be derived by formulas (\ref{2.28}) and (\ref{2.29}):
\begin{equation}
\bR_X(\omega)=\bT^{-1} \bR_{E_O} (\omega) \bT:~(r_i^j(X,\omega)).
\label{2.30}
\end{equation}

After \cite{M97}, we introduce the so-called hyperboloid parametrization as follows
\begin{equation}
\begin{gathered}
x^0=\cosh{r} \cos{\phi}, ~ ~
x^1=\cosh{r} \sin{\phi}, \\
x^2=\sinh{r} \cos{(\theta-\phi)}, ~ ~
x^3=\sinh{r} \sin{(\theta-\phi)},
\end{gathered} \label{2.31}
\end{equation}
where $(r,\theta)$ are the polar coordinates of the base plane and $\phi$ is just the fibre coordinate. We note that
$$-x^0x^0-x^1x^1+x^2x^2+x^3x^3=-\cosh^2{r}+\sinh^2{r}=-1<0.$$
The inhomogeneous coordinates corresponding to (2.9), that play an important role in the later visualization of prism tilings in $\EUC$,
are given by
\begin{equation}
\begin{gathered}
x=\frac{x^1}{x^0}=\tan{\phi}, ~ ~
y=\frac{x^2}{x^0}=\tanh{r} \frac{\cos{(\theta-\phi)}}{\cos{\phi}}, \\
z=\frac{x^3}{x^0}=\tanh{r} \frac{\sin{(\theta-\phi)}}{\cos{\phi}}.
\end{gathered} \label{2.32}
\end{equation}
\subsubsection{Translation curves}
We recall some basic facts about translation curves in $\SLR$  following 
\cite{MSz, MSzV, MSzV17}. For any point $X(x^{0}, x^{1}, x^{2}, x^{3}) \in {\mathcal H}$ 
(and later also for points in $\widetilde{\mathcal H}$) the \emph{translation map} from  the origin 
$E_{0} (1;0;0;0)$ to $X$ is defined by the \emph{translation matrix} $\bT$ and its inverse presented in (\ref{2.28}).

Let us consider for a given vector  $(q;u;v;w)$ a curve 
$\mathcal C(t) = (x^{0}(t);x^{1}(t);$ $x^{2}(t); x^{3}(t))$, $t \geq 0$, in ${\mathcal H}$ starting at the origin: $\mathcal C(0) = E_{0}(1;0;0;0)$ and such that
$$
\dot{\mathcal C} (0) = (\dot{x}^{0}(0), \dot{x}^{1}(0), \dot{x}^{2}(0), \dot{x}^{3}(0)) = (q,u,v,w),
$$
where $\dot{\mathcal C}(t)  = (\dot{x}^{0}(t), \dot{x}^{1}(t), \dot{x}^{2}(t), \dot{x}^{3}(t))$
is the tangent vector at any point of the curve. For $t \geq 0$ there exists a matrix
\begin{equation*}
\bT(t) = \left( \begin{array}{cccc}
x^{0}(t) & x^{1}(t) & x^{2}(t) & x^{3}(t) \cr
-x^{1}(t) & x^{0}(t) & x^{3}(t) & -x^{2}(t) \cr
x^{2}(t) & x^{3}(t) & x^{0}(t) & x^{1}(t) \cr
x^{3}(t) & -x^{2}(t) & -x^{1}(t) & x^{0}(t) \cr
\end{array} \right) \tag{2.27}
\end{equation*}
which defines the translation from $\mathcal C(0)$ to $\mathcal C (t)$:
\begin{equation*}
\mathcal C(0) \cdot \bT(t) = \mathcal C(t),  \quad t \geqslant 0. \tag{2.28}
\end{equation*}
The $t$-parametrized family $\bT(t)$ of translations is used in the following definition.

As we mentioned this earlier, 
the curve $\mathcal C(t)$, $t \geqslant 0$, is said to be a \emph{translation curve} if
\begin{equation*}
\dot{\mathcal C }(0) \cdot \bT(t) = \dot{\mathcal C}(t),  \quad t \geqslant 0. \tag{2.29}
\end{equation*}

The solution, depending on $(q,u,v,w)$ had been determined in \cite{MoSzi10}, where it splits into three cases. 

It was observed above that for any $X(x^{0},x^{1},x^{2},x^{3}) \in \widetilde{\mathcal H}$ there 
is a suitable transformation $\bT^{-1}$, which sent $X$ to the origin $E_{0}$ along a translation curve.
For a given translation curve $\mathcal C = \mathcal C (t)$ the initial unit tangent vector $(u,v,w)$ (in Euclidean coordinates) 
at $E_{0}$ can be presented as
\begin{equation}
u = \sin \alpha, \quad v = \cos \alpha \cos \lambda, \quad w = \cos \alpha \sin \lambda, \notag
\end{equation}
for some $-\frac{\pi}{2} \leqslant \alpha \leqslant \frac{\pi}{2}$ and $ -\pi < \lambda \leqslant \pi$.  
In $\widetilde{\mathcal H}$ this vector is of length square $-u^{2} + v^{2} + w^{2}  = \cos 2 \alpha$. We 
always can assume that $\mathcal C$ is parametrized by the translation arc-length parameter $t = s \geqslant 0$. 
Then coordinates of a point $X(x; y; z)$ of $\mathcal C$, such that the translation distance between $E_{0}$ and $X$ equals 
$s$, depend on $(\lambda,\alpha,s)$ as geographic coordinates according to the above considered three cases as follows.
\begin{table}[ht]
\caption{Translation  curves.} \label{table2}
\vspace{3mm}
\centerline{
$
\begin{array}{|c|l|} \hline
\textrm{direction} & \textrm{parametrization of a translation curve}  \cr \hline  
\begin{gathered}
0 \le \alpha < \frac{\pi}{4} \\ (\bH^2-{\rm like})
\end{gathered}
&
\begin{array}{l}
\begin{gathered} \noalign{\vskip2pt}\begin{pmatrix} x(s,\alpha,\lambda)\\y(s,\alpha,\lambda)\\z(s,\alpha,\lambda)\end{pmatrix}=  \frac{\tanh (s \sqrt{\cos 2 \alpha})}{\sqrt{\cos 2\alpha}} \begin{pmatrix} \sin \alpha\\
 \cos \alpha \cos \lambda\\
 \cos \alpha \sin \lambda)\end{pmatrix} \end{gathered} \cr \noalign{\vskip2pt}\end{array} \cr
 \hline
\begin{gathered} \alpha=\frac{\pi}{4} \\ ({\rm light-like}) \end{gathered}
&
\begin{array}{l}
\begin{gathered} \noalign{\vskip2pt}\begin{pmatrix} x(s,\alpha,\lambda)\\y(s,\alpha,\lambda)\\z(s,\alpha,\lambda)\end{pmatrix}=  \frac{\sqrt{2} s}{2} 
\begin{pmatrix} 1\\
 \cos \lambda\\
 \sin \lambda)\end{pmatrix} \end{gathered} \cr
\noalign{\vskip2pt}\end{array} \cr \hline
\begin{gathered} \frac{\pi}{4}  < \alpha \le \frac{\pi}{2} \\ ({\rm fibre-like}) \end{gathered}
&
\begin{array}{l}  \begin{gathered} \noalign{\vskip2pt}\begin{pmatrix} x(s,\alpha,\lambda)\\y(s,\alpha,\lambda)\\z(s,\alpha,\lambda)\end{pmatrix}=  \frac{\tan (s \sqrt{-\cos 2 \alpha})}{\sqrt{-\cos 2\alpha}} \begin{pmatrix} \sin \alpha\\
 \cos \alpha \cos \lambda\\
 \cos \alpha \sin \lambda)\end{pmatrix} \end{gathered}   \cr
\noalign{\vskip2pt}
\end{array}  \cr \hline
\end{array}
$
}
\end{table}
\section{Translation-like Apollonius surfaces}
In $\SXR$ and $\HXR$ geometries, the translation and geodesic curves coincide with each other and in the paper \cite{Sz22} we examined and visualized them. Moreover, we generalized the corresponding Ceva and Menelaus theorems in the mentioned spaces (see \cite{Sz22-3}). 
One of our further goals is to examine and visualize the Dirichlet-Voronoi cells of $X \in \{\NIL, \SOL, \SLR\}$ geometry. 
In order to get D-V cells we have to determine its ``faces" that are
parts of bisector (or equidistant) surfaces of given point pairs. We determine the equations of a larger group of surfaces than the translation-like Apollonius surfaces that are defined in Definition 1.1. The translation distances in $X$ geometry are given by the formerly determined translation curves.
\begin{defn} \label{def-distance}
A  \emph{translation distance in $X$ geometry} $\rho (E_{0},P)$ between the origin $E_{0}(1,0,0,0)$ and the point  $P(1,a,b,c)$ is the length of a 
translation curve connecting them.
\end{defn}
In each investigated geometry, the key to solving the question is the so-called ``inverse problem", in which the appropriate parameters of the translation curve starting from the starting point are found for a known point $P=(1,a,b,c)$. 
\subsection{$\NIL$ geometry}
\subsubsection{Inverse problem}
It can be assumed by the homogeneity of $\NIL$ that the starting point of a 
given translation curve segment is $E_0=P_1=(1,0,0,0)$ and 
the other endpoint will be given by its homogeneous coordinates $P_2=(1,a,b,c)$. 
We consider the translation curve segment $t_{P_1P_2}$ and determine its
parameters $(\phi,\theta,r)$ expressed by the real coordinates $a$, $b$, $c$ of $P_2$. 
We obtain directly by equation system (\ref{2.14}) the following:
\begin{lem}
\begin{enumerate}
\item Let $(1,a,b,c)$ $(a,b \in \mathbf{R} \setminus \{0\},~c\in \mathbf{R})$ be the homogeneous 
coordinates of the point $P \in \NIL$. The parameters of the
corresponding translation curve $t_{E_0P}$ are the following
\begin{equation}
\begin{gathered}
\phi=\mathrm{arccot}\Big(\frac{a}{b}\Big),~\text{or}~ \phi=\mathrm{arccot}\Big(\frac{a}{b}\Big)-\pi, \\
\theta=\mathrm{arctan}\Big( \frac{c-\frac{ab}{2}}{\sqrt{a^{2}+b^{2}}}\Big),~
r=\Big|\frac{c-\frac{ab}{2}}{\sin{\theta}}\Big|.
\end{gathered} \label{3.1}
\end{equation}
\item Let $(1,a,0,c)$ $(a,c \in \mathbf{R} \setminus \{0\})$ be the homogeneous 
coordinates of the point $P \in \NIL$. The parameters of the
corresponding translation curve $t_{E_0P}$ are the following
\begin{equation}
\begin{gathered}
\phi=\pi \cdot n ,~ (n\in\{0,1\}),~ 
\theta= \mathrm{arctan}\Big(\frac{c}{a}\Big),~
r=\Big|\frac{a}{\cos{\theta}}\Big|.
\end{gathered} \label{3.2}
\end{equation}
\item Let $(1,a,0,0)$ $(a \in \mathbf{R}\setminus \{0\})$ 
be the homogeneous coordinates of the point $P \in \NIL$. The parameters of the 
corresponding translation curve $t_{E_0P}$ are the following
\begin{equation}
\begin{gathered}
\phi=\pi \cdot n ,~ (n\in\{0,1\}),~
\theta=\pi \cdot n ,~ (n\in\{0,1\}),~
r=|a|.
\end{gathered} \label{3.3}
\end{equation}
\item Let $(1,0,b,0)$ $(b \in \mathbf{R}\setminus \{0\})$ 
be the homogeneous coordinates of the point $P \in \NIL$. The parameters of the 
corresponding translation curve $t_{E_0P}$ are the following
\begin{equation}
\begin{gathered}
\phi=\pm \frac{\pi}{2},~
\theta=\pi \cdot n ,~ (n\in\{0,1\}),~
r=|b|.
\end{gathered} \label{3.4}
\end{equation}
\item Let $(1,0,0,c)$ $(c \in \mathbf{R}\setminus \{0\})$ 
be the homogeneous coordinates of the point $P \in \NIL$. The parameters of the 
corresponding translation curve $t_{E_0P}$ are the following
\begin{equation}
\begin{gathered}
\theta=\pm \frac{\pi}{2},~
r=|c|.~ ~ \square
\end{gathered} \label{3.5}
\end{equation}
\end{enumerate}
\label{lem1}
\end{lem}
\subsection{$\SOL$ geometry}
\subsubsection{Inverse problem}
It can be assumed by the homogeneity of $X$ that the starting point of a given translation curve segment is $E_0=P_1=(1,0,0,0)$.
The other endpoint will be given by its homogeneous coordinates $P_2=(1,a,b,c)$. We consider the translation curve segment $t_{P_1P_2}$ and determine its
parameters $(\phi,\theta,t)$ expressed by the real coordinates $a$, $b$, $c$ of $P_2$. We obtain directly by equation system (\ref{2.26}) the following Lemma (see \cite{Sz17} and \cite{Sz19}):
\begin{lem}
\begin{enumerate}
\item Let $(1,a,b,c)$ $(b,c \in \bR \setminus \{0\}, a \in \bR)$ be the homogeneous coordinates of the point $P \in \SOL$. The parameters of the
corresponding translation curve $t_{E_0P}$ are the following
\begin{equation}
\begin{gathered}
\phi=\mathrm{arccot}\Big(-\frac{a}{b} \frac{\mathrm{e}^{c}-1}{\mathrm{e}^{-c}-1}\Big),~\theta=\mathrm{arccot}\Big( \frac{b}{\sin\phi(\mathrm{e}^{c}-1)}\Big),\\
t=\frac{c}{\sin\theta}, ~ \text{where} ~ -\pi < \phi \le \pi, ~ -\pi/2\le \theta \le \pi/2, ~ t\in \bR^+.
\end{gathered} \label{3.6}
\end{equation}
\item Let $(1,a,0,c)$ $(a,c \in \bR \setminus \{0\})$ be the homogeneous coordinates of the point $P \in \SOL$. The parameters of the
corresponding translation curve $t_{E_0P}$ are the following
\begin{equation}
\begin{gathered}
\phi=0~\text{or}~  \pi, ~\theta=\mathrm{arccot}\Big( \mp \frac{a}{(\mathrm{e}^{-c}-1)}\Big),\\
t=\frac{c}{\sin\theta}, ~ \text{where}  ~ -\pi/2\le \theta \le \pi/2, ~ t\in \bR^+.
\end{gathered} \label{3.7}
\end{equation}
\item Let $(1,a,b,0)$ $(a,b \in \bR)$ be the homogeneous coordinates of the point $P \in \SOL$. The parameters of the 
corresponding translation curve $t_{E_0P}$ are the following
\begin{equation}
\begin{gathered}
\phi=\arccos\Big(\frac{x}{\sqrt{a^2+b^2}}\Big),~  \theta=0,\\
t=\sqrt{a^2+b^2}, ~ \text{where}  ~ -\pi < \phi \le \pi, ~ t\in \bR^+.~ ~ \square
\end{gathered} \label{3.8}
\end{equation}
\end{enumerate}
\label{invlemsol}
\end{lem}
\subsection{$\SLR$ geometry}
\subsubsection{Inverse problem}
It can be assumed by the homogeneity of $\SLR$ that the starting point of a 
given translation curve segment is $E_0=P_1=(1,0,0,0)$ and 
the other endpoint will be given by its homogeneous coordinates $P_2=(1,a,b,c)$. 
We consider the translation curve segment $t_{P_1P_2}$ and determine its
parameters $\alpha,\lambda,s$ expressed by the real coordinates $a$, $b$, $c$ of $P_2$.We can assume by the structure of $\SLR$ that $a \in \bR^+$.
We obtain directly by equation systems in Table 1. the following:
\begin{lem}
\begin{enumerate}
\item Let $(1,a,b,c)$ $(0 \ne b,c \in \bR, a \in \bR, a^2-b^2-c^2<0\})$ $(0 \le \alpha < \frac{\pi}{4}, ~ -\pi < \lambda \le \pi)$ be the homogeneous coordinates of the point $P \in \SLR$. The parameters of the
corresponding translation curve $t_{E_0P}$ are the following
\begin{equation}
\begin{gathered}
\alpha=\mathrm{arctg}\Big(\frac{1}{\sqrt{b^2+c^2}}|a| \Big), \ \
\lambda=\mathrm{arctg}\Big(\frac{c}{b}\Big), \\
s=\frac{\mathrm{arctanh}\Big(\sqrt{a^2+b^2+c^2}\cdot \sqrt{\cos\Big(2 \cdot \mathrm{arctan} \Big(\frac{a}{\sqrt{b^2+c^2}}\Big)\Big)}}{\sqrt{\cos\Big(2\mathrm{arctan} \Big(\frac{a}{\sqrt{b^2+c^2}}\Big) \Big)}}
\end{gathered} \label{3.9}
\end{equation}
\item Let $(1,a,b,c)$ $(0 \ne b,c \in \bR, a \in \bR^+, a^2-b^2-c^2=0\})$ $-\pi < \lambda \le \pi)$ be the homogeneous coordinates of the point $P \in \SLR$. The parameters of the
corresponding translation curve $t_{E_0P}$ are the following
\begin{equation}
\begin{gathered}
\alpha=\frac{\pi}{4}, \ \ \lambda=\mathrm{arctg}\Big(\frac{c}{b}\Big), \ \ \ s=\frac{a}{\sqrt{2}} ,\\
\end{gathered} \label{3.10}
\end{equation}
\item Let $(1,a,b,c)$ $(0 \ne b,c \in \bR, a \in \bR^+, a^2-b^2-c^2>0\})$ $(\frac{\pi}{4} < \alpha < \frac{\pi}{2}, ~ -\pi < \lambda \le \pi)$ be the homogeneous coordinates of the point $P \in \SLR$. The parameters of the
corresponding translation curve $t_{E_0P}$ are the following
\begin{equation}
\begin{gathered}
\alpha=\mathrm{arctg}\Big(\frac{1}{\sqrt{b^2+c^2}}a \Big), \ \
\lambda=\mathrm{arctg}\Big(\frac{c}{b}\Big), \\
s=\frac{\mathrm{arctanh}\Big(\sqrt{a^2+b^2+c^2}\cdot \sqrt{-\cos\Big(2 \cdot \mathrm{arctan} \Big(\frac{a}{\sqrt{b^2+c^2}}\Big)\Big)}}{\sqrt{-\cos\Big(2\mathrm{arctan} \Big(\frac{a}{\sqrt{b^2+c^2}}\Big) \Big)}}
\end{gathered} \label{3.11}
\end{equation}
\item Let $(1,a,b,c)$ $(0 = b,c \in \bR, a \in \bR^+\})$ be the homogeneous coordinates of the point $P \in \SLR$. The parameters of the
corresponding translation curve $t_{E_0P}$ are the following
\begin{equation}
\begin{gathered}
\alpha=\frac{\pi}{2}, \ \ s=\mathrm{arctg}(x).
\end{gathered} \ \ \ \ \square \label{3.12}
\end{equation}
\end{enumerate}
\label{invlem}
\end{lem}
\subsection{The equations of Apollonius surfaces}
The Apollonius surfaces are determined in a similar way in all three investigated geometries. The basic principle, which is based on the application of the isometries and the corresponding ``inverse problem" of the given space, is illustrated in the $\NIL$ geometry.

{\it In order to determine the translation-like Apollonius surface
$\cA\cS_{P_1P_2}$ related to a translation segment $t_{P_1P_2}$ . Let P=(1,x,y,z) be the point lying on the Apollonius surface. we define {translations} $\bT_{P_2}$, as elements of the isometry group of $\NIL$, that
maps the origin $E_0$ onto $P$ where $P_1=E_0=(1,0,0,0)$ and $P_2=(1,a,b,c)$}.

The isometry $\bT_{P_2}$ and its inverse (up to a positive determinant factor) can be given by:
\begin{equation}
\bT_{P_2}=
\begin{pmatrix}
1 &a & b & c \\
0 & 1 & 0 & 0 \\
0 & 0 & 1 & a \\
0 & 0 & 0 & 1
\end{pmatrix} , ~ ~ ~
\bT_{P_2}^{-1}=
\begin{pmatrix}
1 & -a & -b & ab - c \\
0 & 1 & 0 & 0 \\
0 & 0 & 1 & -a \\
0 & 0 & 0 & 1 
\end{pmatrix}, \notag
\end{equation}
and the image $P^{2}=\bT^{-1}_{P_2}(P)$ of the point $P$ is the following :
\begin{equation}
\begin{gathered}
\bT^{-1}_{P_2}(P)=(1,-a+x,-b+y,a(b-y)-c+z); \\ \text{and it is obvious that } \bT^{-1}_{P_2}(P_2)=E_0=(1,0,0,0).
\end{gathered} \label{3.13}
\end{equation}
It is clear that $P=(1,x,y,z) \in \cA\cS^\lambda_{P_1P_2} ~ \text{iff} ~ \sigma \cdot d^t(P_1,P)=  d^t(P,P_2) \Rightarrow \sigma \cdot d^t(P_1,P)= \cdot d^t(P_1=E_0,P^2)$ where
$P^2=\bT^{-1}_{P_2}(P)$ and $\sigma \in\bR^+$.

This method leads to equations of the Apollonius surfaces in investigated spaces.
\begin{lem}
The implicit equation of the translation-like Apollonius surface $\cA\cS^\sigma_{P_1P_2}(x,y,z)$ of two points $P_1=(1,0,0,0)$, $P_2=(1,a,b,c)$ and with parameter $\sigma \in\bR^+$ in $\NIL$ space (see Fig.~1-2.): 
\begin{equation}
\begin{gathered}
\sigma \cdot \sqrt{((x^2+y^2)+(\frac{1}{2}xy-z)^2)}- \\ 
\sqrt{((-x+a)^2+(-y+b)^2+ (-xb+xy+c-z-\frac{1}{2}(-x+a) \cdot (-y+b))^2}=0.
\end{gathered} \label{3.14}
\end{equation}
\end{lem}
\begin{figure}[ht]
\centering
\includegraphics[width=0.98\linewidth]{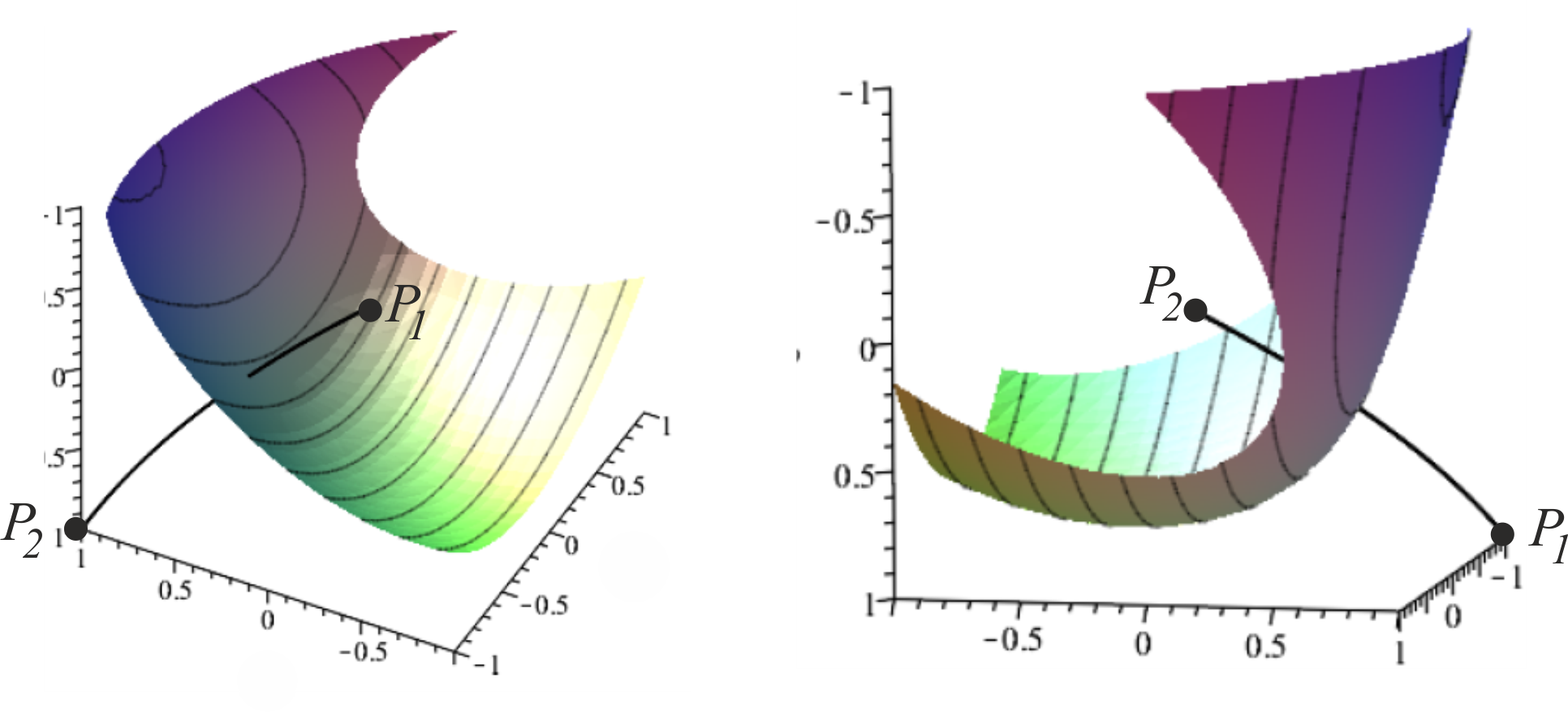}
\caption{Translation-like $\NIL$ Apollonius surface of point pairs $(P_1,P_2)$ with coordinates $((1,0,0,0), (1,-1,1,1))$ with parameter $\lambda=2$,}
\label{pic:surf4}
\end{figure}
\begin{figure}[ht]
\centering
\includegraphics[width=0.98\linewidth]{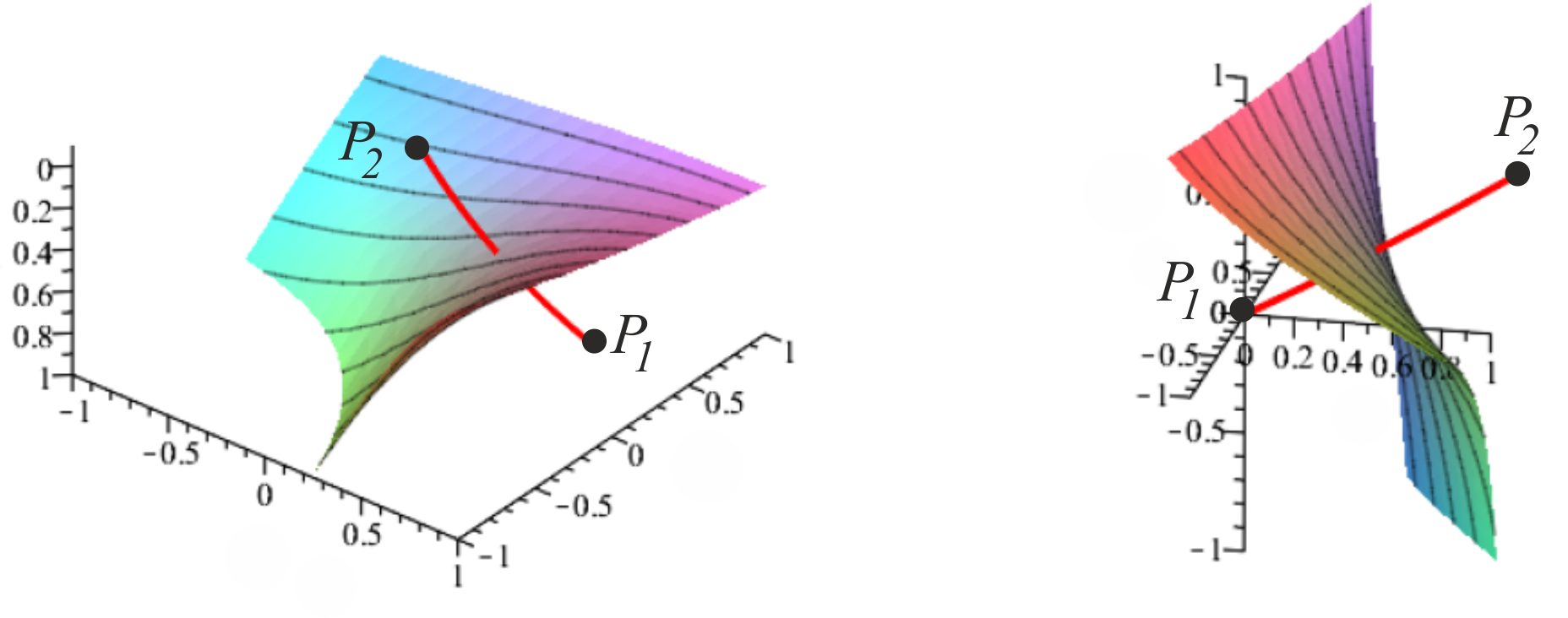}
\caption{Translation-like $\NIL$ bisector (equidistant surface) of point pairs $(P_1,P_2)$ with coordinates $((1,0,0,0), (1,1/2,1,1/2))$. }
\label{pic:surf3}
\end{figure}
\begin{lem}
The implicit equation of the translation-like Apollonius surface $\cA\cS^\sigma_{P_1P_2}(x,y,z)$ of two points $P_1=(1,0,0,0)$, $P_2=(1,a,b,c)$ and with parameter $\lambda \in\bR^+$ in $\SOL$ space (see Fig.~3-4.):
\begin{enumerate}
\item $c \ne 0$
\begin{equation}\label{bis1}
\begin{gathered}
z\ne 0, c~:~
\frac{|c-z|}{|\mathrm{e}^{c}-\mathrm{e}^z|}\sqrt{(a-x)^2 \mathrm{e}^{2(c+z)}+(\mathrm{e}^{c}-\mathrm{e}^{z})^2+(b-y)^2}=\\
=\sigma \cdot \frac{|z|}{|\mathrm{e}^{z}-1|}\sqrt{x^2 \mathrm{e}^{2z}+(\mathrm{e}^{z}-1)^2+y^2},\\
z=c~:~\sqrt{(x-a)^2\mathrm{e}^{2c}+(y-b)^2\mathrm{e}^{-2c}}
=\sigma \cdot \frac{|z|}{|\mathrm{e}^{z}-1|}\sqrt{x^2 \mathrm{e}^{2z}+(\mathrm{e}^{z}-1)^2+y^2},\\
z=0~:~ \frac{|c|}{|\mathrm{e}^{c}-1|}\sqrt{(a-x)^2 \mathrm{e}^{2c}+(\mathrm{e}^{c}-1)^2+(b-y)^2}=\lambda \cdot \sqrt{x^2+y^2},
\end{gathered} \tag{3.6}
\end{equation}
\item $c=0$
\begin{equation}\label{bis1a}
\begin{gathered}
z\ne 0~:~
\frac{|z|}{|\mathrm{e}^z-1|}\sqrt{(a-x)^2 \mathrm{e}^{2z}+(\mathrm{e}^{z}-1)^2+(b-y)^2}=\\
=\sigma \cdot \frac{|z|}{|\mathrm{e}^{z}-1|}\sqrt{x^2 \mathrm{e}^{2z}+(\mathrm{e}^{z}-1)^2+y^2},\\
z=0~:~ \sqrt{(x-a)^2+(y-b)^2}=\sigma \cdot \sqrt{x^2+y^2}.~\square
\end{gathered} \tag{3.7}
\end{equation}
\end{enumerate}
\end{lem}
\begin{figure}[ht]
\centering
\includegraphics[width=0.98\linewidth]{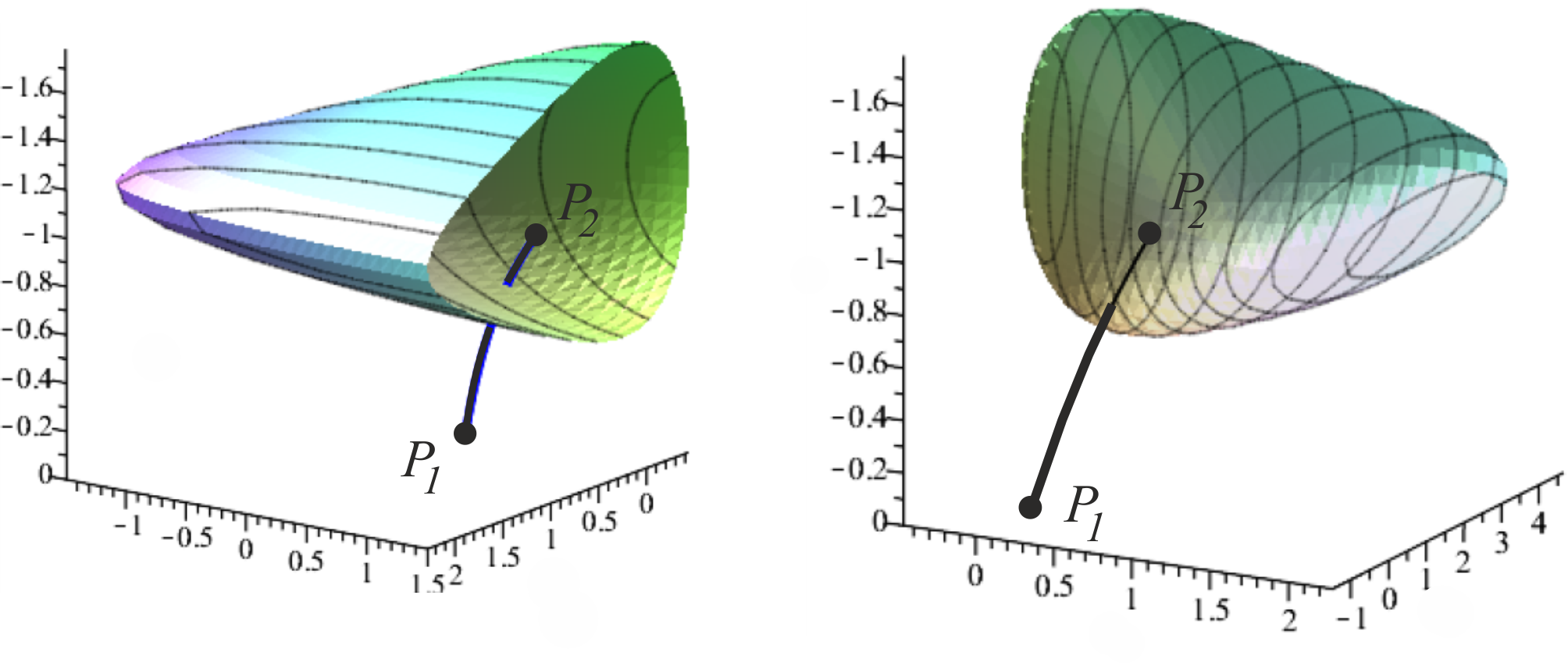}
\caption{Translation-like$ \SOL$ Apollonius surface of point pairs $(P_1,P_2)$ with coordinates $((1,0,0,0), (1,-1,1,1/2))$ with parameter $\sigma=1/2$,}
\label{pic:surf1}
\end{figure}
\begin{figure}[ht]
\centering
\includegraphics[width=0.98\linewidth]{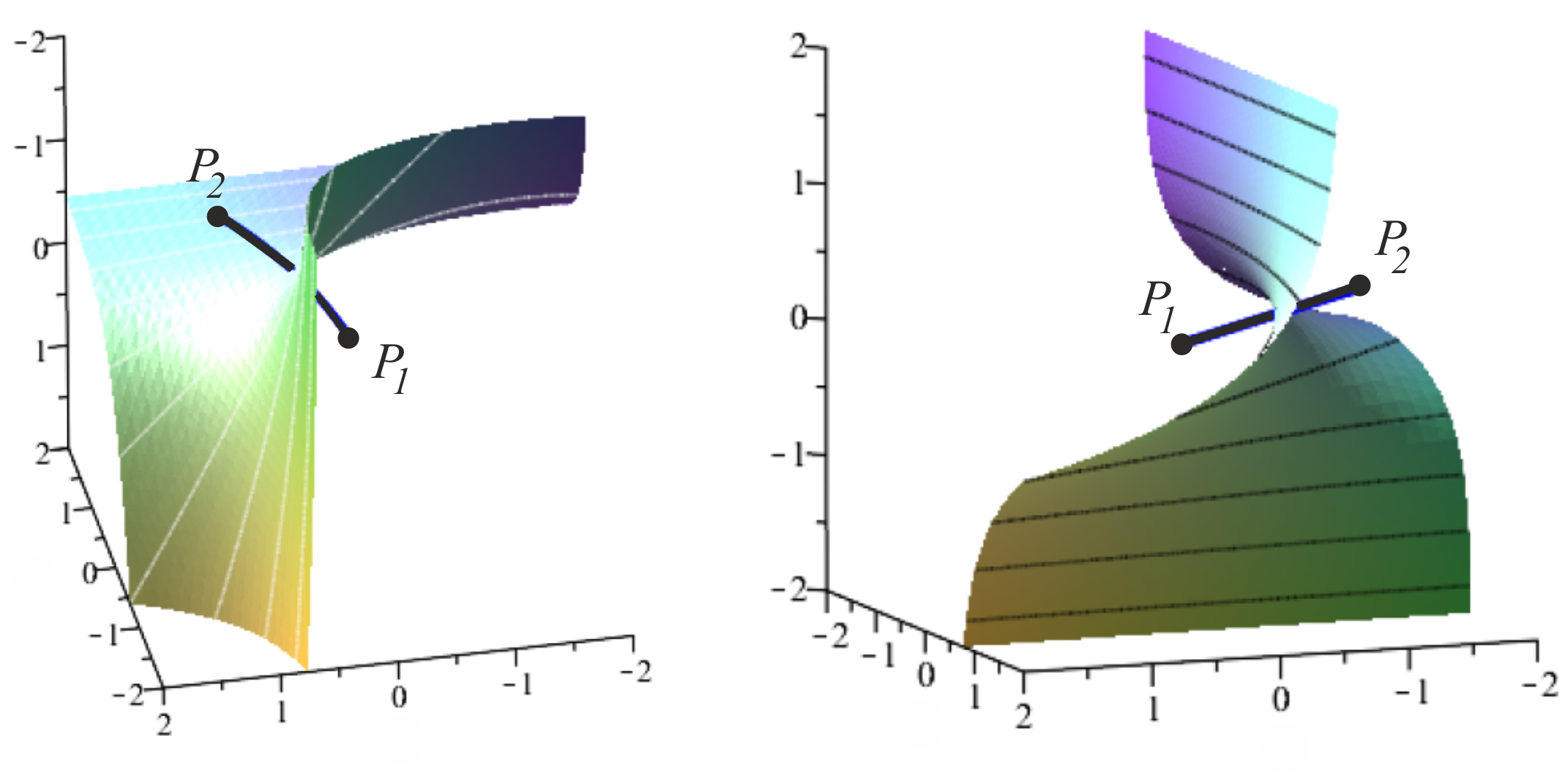}
\caption{Translation-like $\SOL$ bisector (equidistant surface) of point pairs $(P_1,P_2)$ with coordinates $((1,0,0,0), (1,-1,1,1/2))$. }
\label{pic:surf2}
\end{figure}
\begin{lem}
The implicit equation of the translation-like Apollonius surface $\cA\cS^\sigma_{P_1P_2}(x,y,z)$ of two points $P_1=(1,0,0,0)$, $P_2=(1,a,b,c)$ and with parameter $\sigma \in\bR^+$ in $\SLR$ space is $\sigma \cdot d^t(P_1,P)= d^t(P_1=E_0,P^2)$. 
(Unfortunately, this equation is very long and we will not publish it here, but it can be easily reproduced based on the examples of the previous spaces and Lemma 3.4 and moreover this is described in Fig.~5-6).
\end{lem}
\begin{figure}[ht]
\centering
\includegraphics[width=0.98\linewidth]{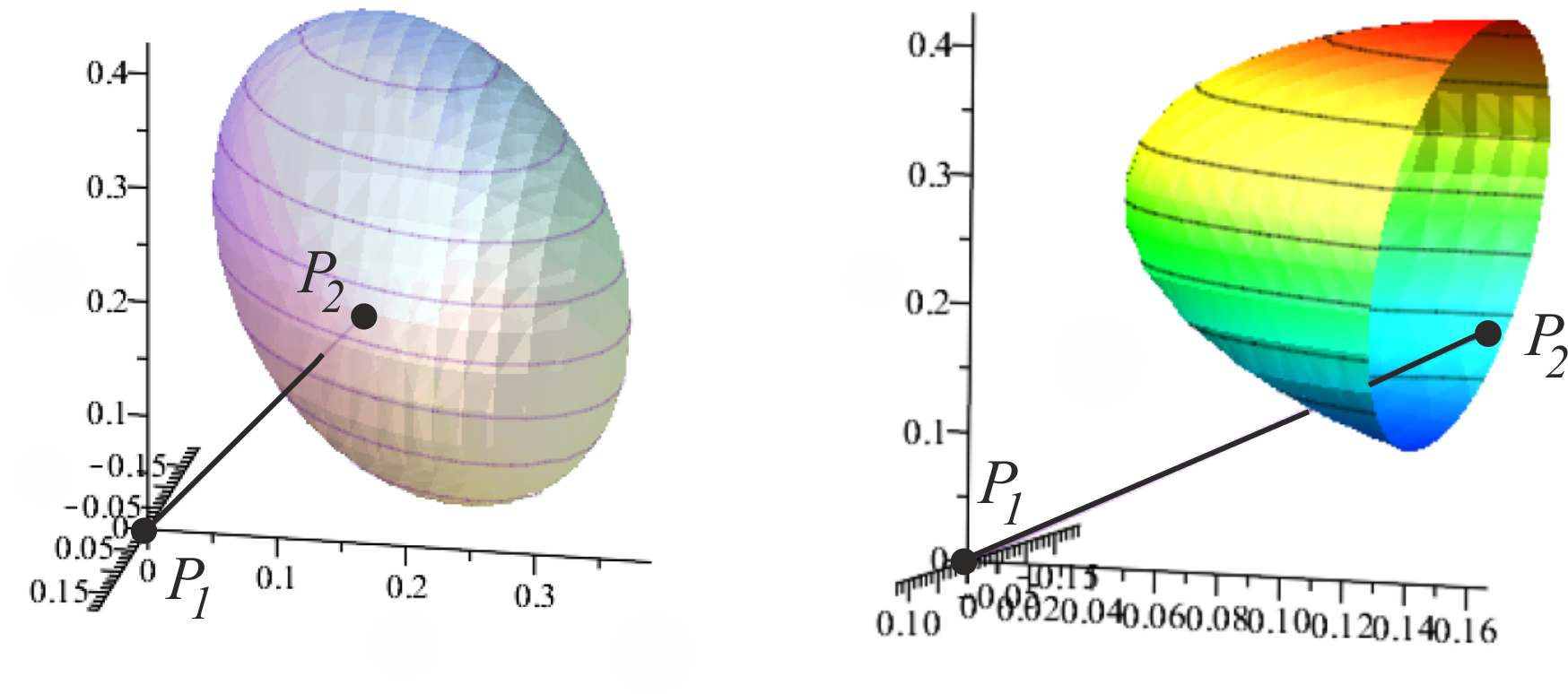}
\caption{Translation-like $\SLR$ Apollonius surface of point pairs $(P_1,P_2)$ with coordinates $((1,0,0,0), (1,0,1/6,1/5))$ with parameter $\sigma=2$,}
\label{6}
\end{figure}
\begin{figure}[ht]
\centering
\includegraphics[width=0.98\linewidth]{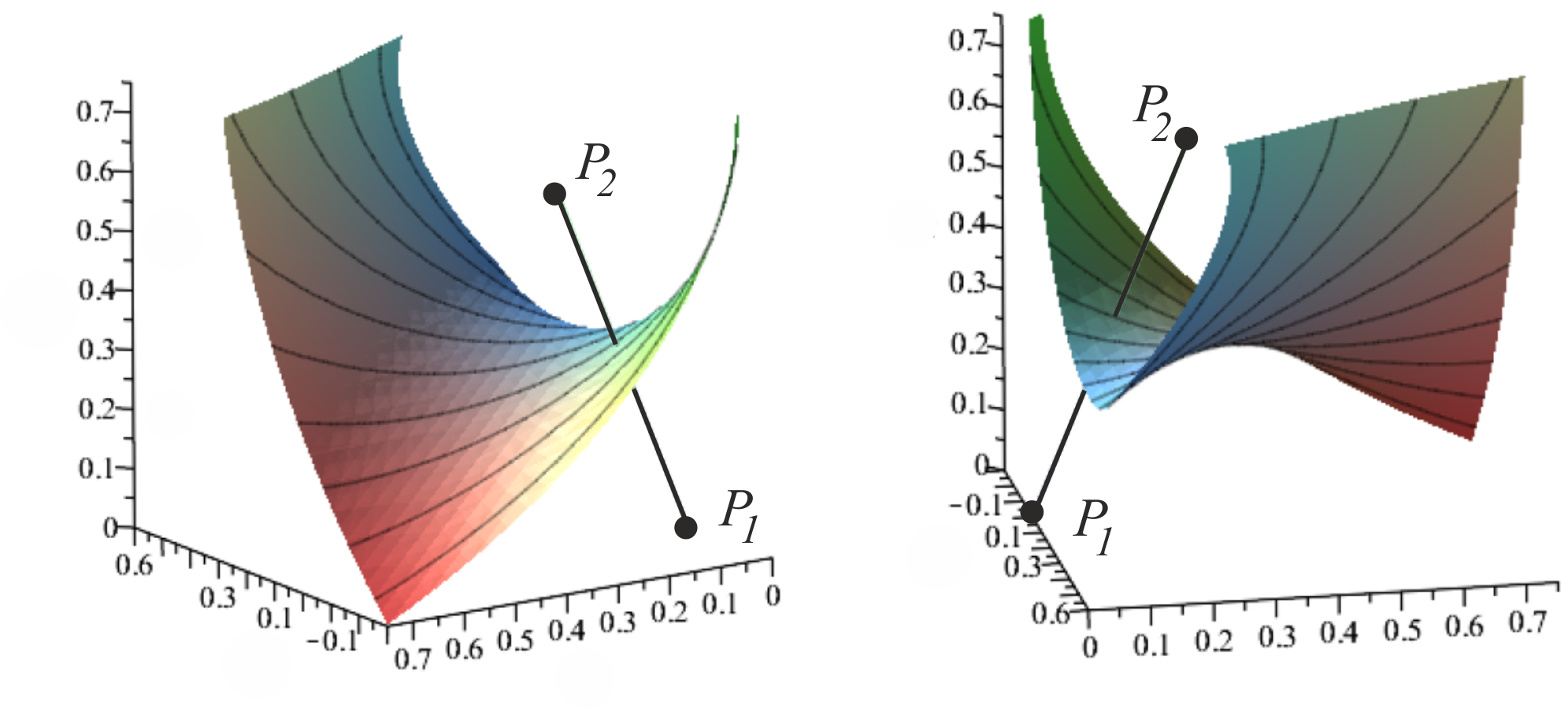}
\caption{Translation-like $\SLR$ bisector (equidistant surface) of point pairs $(P_1,P_2)$ with coordinates $((1,0,0,0), (1,1/4,3/5,0))$. }
\label{pic:surf5}
\end{figure}
\section{Surface of translation-like triangle}
In the previous sections, we gave the ``plane"-like surface in such a way that we tried to exploit the generalized geometric locus problem related to the Apollonius (or in a special case, equidistant) surface of two given points. In this section, we approach the notion of the generalized plane from another possible direction. According to an axiom in Euclidean geometry, three non-collinear points $P_1$, $P_2$ and $P_3$ determine uniquely a plane. Furthermore, this is \textit{transitive}, i.e. any three different points $P_1'$, $P_2'$ and $P_3'$ of this plane define the same plane. There are countless ways to check whether a given point $P$ is located on the plane or not. We choose from these whether the lines $PP_1$, $PP_2$ and $PP_3$ are coplanar or not. Of course, in the given $X$ Thurston geometry, these can also be translation or geodesic curves. We can also speak of coplanarity only locally, at the point $P$. Based on these, we can formulate the following definition:

\begin{defn}
    The $S^{X,t}_{P_1,P_2,P_3}$ translation-like triangular surface of $P_1,$ $P_2$ and $P_3$ in the Thurston geometry $X \in \{\EUC,\SPH,\HYP,\SXR,\HXR,\NIL,\SLR,\SOL\}$ is  
the set of all points of $X$ from which the tangents in $P$ of the translation curves drawn to points $P_1,$ $P_2$ and $P_3$ are coplanar.
\label{def:plane}
\end{defn}

If we consider the Thurston geometries where the translation curves are represented as Euclidean lines in the projective model of the corresponding geometry $(\EUC,\ \HYP,\ \SPH,\ \SLR )$, then the set of points obtained by Definition \ref{def:plane} is uniquely represented as a Euclidean plane. It is also easy to see that the transitivity property is automatically satisfied in these geometries. In the following, we will review geometries that do not have these properties.

\subsection{Translation-like triangular surface in $\SOL$}

The Apollonius surfaces were examined more thoroughly in $\NIL$ geometry, and the calculations related to translation-like triangular surfaces are presented in detail in $\SOL.$ Let $P_2=(1,a,b,c)$ and $P_3=(1,d,e,f)$  be points in the projective model of $\SOL$ geometry. Due to homogeneity we can assume that $P_1=(1,0,0,0)$ is the origin. At first we need to find the tangents of the translation curves, drawn from point $P=(1,x,y,z)$ to $P_1,$ $P_2$ and $P_3.$ For this, we need to pull back $P$ to the origin. 
     \begin{equation}
     \bT_P=\left(
\begin{array}{cccc}
 1 & x & y & z \\
 0 & e^{-z} & 0 & 0 \\
 0 & 0 & e^z & 0 \\
 0 & 0 & 0 & 1 \\
\end{array}
\right), ~ ~ ~
\bT_{P}^{-1}=
\left(
\begin{array}{cccc}
 1 & -x\,e^z & -y\, e^{-z} & -z \\
 0 & e^z & 0 & 0 \\
 0 & 0 & e^{-z} & 0 \\
 0 & 0 & 0 & 1 \\
\end{array}
\right)
     \end{equation}

\begin{equation}
\begin{gathered}
P^1:=\bT^{-1}_{P}(P_1)=(1,-x\mathrm{e}^{z},-y\mathrm{e}^{-z},-z)\\
P^2:=\bT^{-1}_{P}(P_2)=(1,(a-x)\mathrm{e}^{z},(b-y)\mathrm{e}^{-z},(c-z))\\
P^3:=\bT^{-1}_{P}(P_3)=(1,(d-x)\mathrm{e}^{z},(e-y)\mathrm{e}^{-z},(f-z))    
\end{gathered}
\label{betoltak}
\end{equation}

Now, we need to determine the tangents of the translation curves to $P^i$ $(i=1,2,3)$ at the origin. The most obvious method for this would be to use Lemma \ref{invlemsol}. However, at this point we can simplify the calculations, since we need to derive any nonzero multiple of (\ref{2.25}) from (\ref{2.26}). One can easily see, that the following formula will give  us the unit tangent vector at the origin, multiplied by the $d^t(E_0,P)$ translation distance of $P(1,x,y,z)$ to the origin:

\begin{equation}
   \bt(x,y,z)_\SOL:=d^t(E_0,P)\cdot\bt_e(x,y,z)=\left(\dfrac{x z}{1 - \mathrm{e}^{-z}},\dfrac{y z}{\mathrm{e}^{z}-1}, z\right)
    \label{erintov}
\end{equation}

Applying (\ref{erintov}) to every translated point (\ref{betoltak}), we obtain $\bt(P_i)$ tangent vectors that are coplanar in Euclidean sense iff their triple product is 0.  
\begin{lem}
    Let $P_1=(1,0,0,0),$ $P_2=(1,a,b,c)$ and $P_3=(1,d,e,f)$  be points in the projective model of $\SOL$ geometry. Then the $S^{\SOL,t}_{P_1,P_2,P_3}$ translation-like triangular surface of $P_1,$ $P_2$ and $P_3$ has the equation:
\begin{equation}
    \begin{gathered}
       0=\bt(P^1)\cdot\bt(P^2)\cdot\bt(P^3)= \dfrac{z (z-c) (z-f))}{(1 - \mathrm{e}^{-z}) (\mathrm{e}^z - \mathrm{e}^c) (\mathrm{e}^z -   \mathrm{e}^f)}\cdot\\
       \cdot \left((\mathrm{e}^z-1)(b\, d\, \mathrm{e}^f-a\, e\, \mathrm{e}^c)+   (\mathrm{e}^f-1) (a\, y\,\mathrm{e}^c -b\, x\, \mathrm{e}^z)
        + (\mathrm{e}^c-1) (e\, x\,\mathrm{e}^z -d\, y\,\mathrm{e}^f)\right). 
    \end{gathered}
    \label{solsik}
\end{equation}
\end{lem}
\begin{figure}[h]
    \centering
    \includegraphics[width=0.48\linewidth]{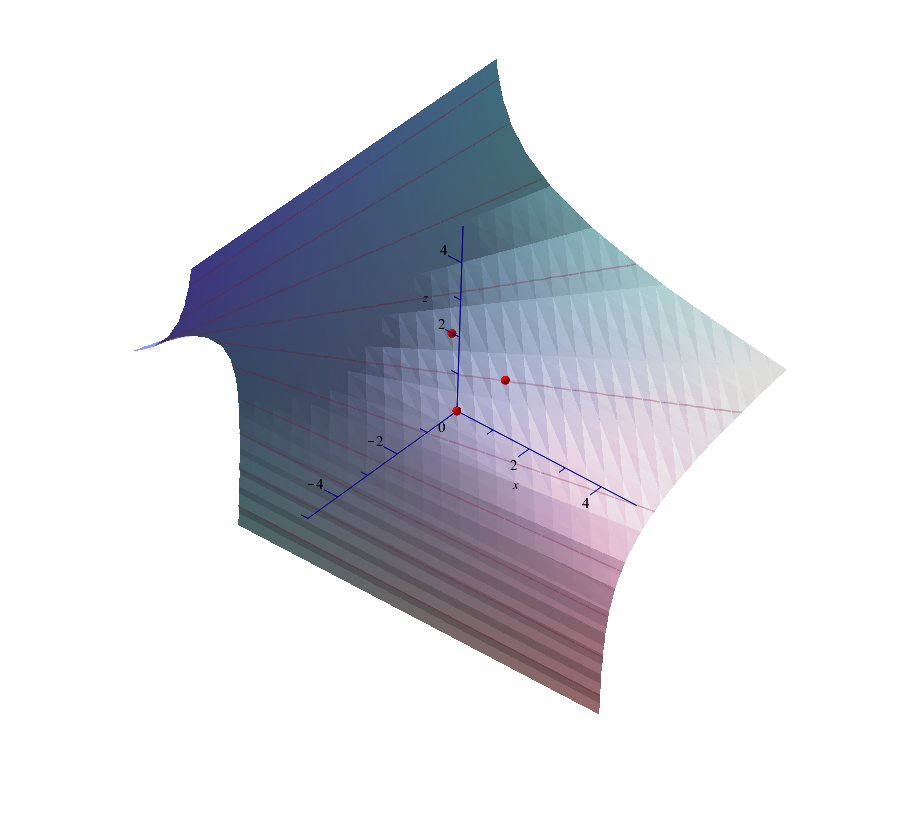}\includegraphics[width=0.48\linewidth]{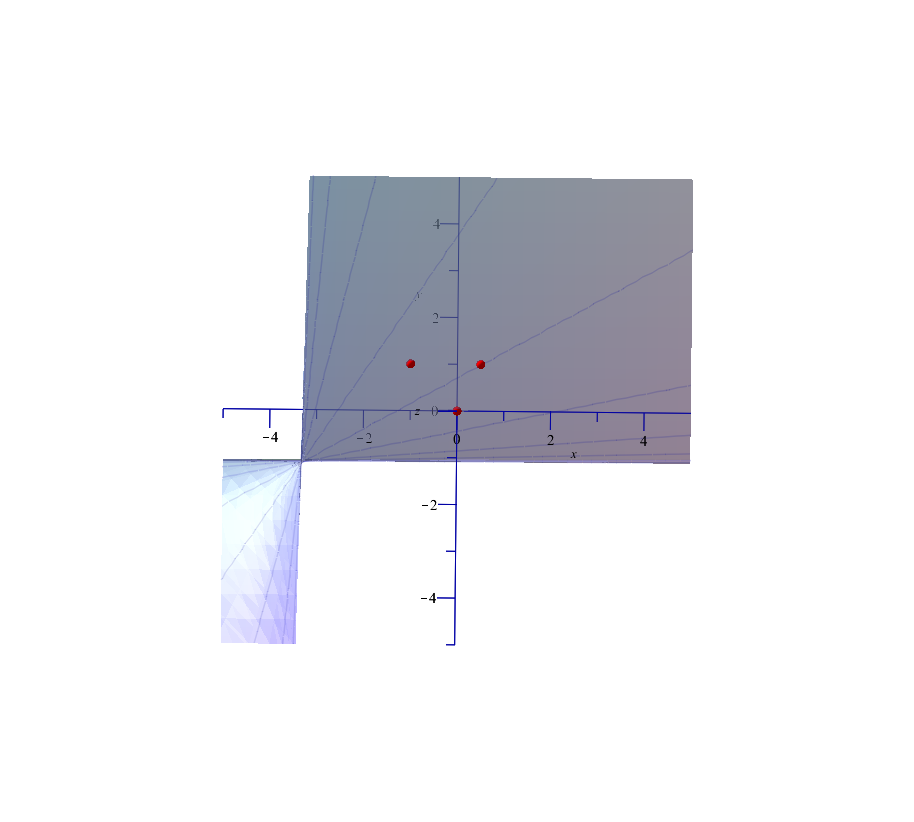}
    \caption{Translation-like triangular surface for $P_1=(1,0,0,0),$ $P_2=(1,-1,1,1)$ and $P_3=(1,\frac{1}{2},1,\frac{1}{2})$}
    \label{fig:Sol_plane}
\end{figure}

An example for the resulted surface can be seen on Figure \ref{fig:Sol_plane}. 
\begin{rem} 
\begin{enumerate}
    \item We can simplify the expression with the fractional multiplier on the right side of (\ref{solsik}), since the limit exists and is finite in all three critical points $(z=0,c,f).$
    \item The level set of (\ref{solsik}) $(z=\mathrm{const})$ are straight lines, \textit{i. e.} the result is a ruled surface. Furthermore the
    $$x=\dfrac{a\, e\, \mathrm{e}^c-b\, d\, \mathrm{e}^f}{e\,(\mathrm{e}^c-1)-b\,(\mathrm{e}^f-1)},\ \ \ 
    y=\dfrac{a\, e\, \mathrm{e}^c-b\, d\, \mathrm{e}^f}
    {d\,\mathrm{e}^f\,(\mathrm{e}^c-1)-a\, \mathrm{e}^c\,(\mathrm{e}^f-1)}$$ 
    line belongs completely to the surface.
\end{enumerate}
\label{solrem}
\end{rem}
\begin{lem}
    Let $P_4$ be a point on $S^{\SOL,t}_{P_1,P_2,P_3}$ translation-like triangular surface. Then $S^{\SOL,t}_{P_1,P_2,P_3}=S^{\SOL,t}_{P_1,P_2,P_4}$ i. e. the definition of the translation-like triangular surface is transitive.
\end{lem}
\textsc{Proof:} It is easy to realise that (\ref{solsik}) equation can be reformulate as an explicit expression for $z$.
This means that the translation-like triangular surface can be considered as a function with 6 parameters and 2 variables:
\begin{equation}
    t^\SOL_{a,b,c,d,e,f}(x,y):=\log{\left(\dfrac{\left( d\,\mathrm{e}^f\,(\mathrm{e}^c-1)-a\, \mathrm{e}^c\,(\mathrm{e}^f-1)\right)\,y-\left(a\, e\, \mathrm{e}^c-b\, d\, \mathrm{e}^f\right)}{\left(e\,(\mathrm{e}^c-1)-b\,(\mathrm{e}^f-1)\right)\, x-\left(a\, e\, \mathrm{e}^c-b\, d\, \mathrm{e}^f\right)}\right)}
    \label{solsikexp}
\end{equation}
Any $P_4$ point on $S^{\SOL,t}_{P_1,P_2,P_3}$ can be considered then as $P_4=\left(1,g,h,t^\SOL_{a,b,c,d,e,f}(g,h)\right)$. We can replace $d,$ $e$ and $f$ in (\ref{solsikexp}) with $g,$ $h$ and $t^\SOL_{a,b,c,d,e,f}(g,h)$ to verify the $t^\SOL_{a,b,c,d,e,f}(x,y)=t^\SOL_{a,b,c,g,h,t^\SOL_{a,b,c,d,e,f}(g,h)}(x,y)$ equation so that this Lemma. This calculation is relatively simple, especially for a computer, but is too extensive to be detailed here. $\ \ \square$

\begin{rem}
    We must pay attention to the domain of $t^\SOL_{a,b,c,d,e,f}(x,y)$. The numerator and denominator of the fraction in the logarithm must have the same sign and none of them can be zero. Fortunately, we have already determined the critical values of $x$ and $y$ in the second paragraph of  Remark \ref{solrem}:
    \begin{equation}
    \begin{gathered}
         \mathrm{dom}(t^\SOL_{a,b,c,d,e,f}(x,y))=\left\{(x,y)\in \mathbb{R}^2\left|
        \left(x-\dfrac{a\, e\, \mathrm{e}^c-b\, d\, \mathrm{e}^f}{e\,(\mathrm{e}^c-1)-b\,(\mathrm{e}^f-1)}\right)\right.\right.\cdot\\
        \cdot
        \left.\left(y-\dfrac{a\, e\, \mathrm{e}^c-b\, d\, \mathrm{e}^f}{d\,\mathrm{e}^f\,(\mathrm{e}^c-1)-a\, \mathrm{e}^c\,(\mathrm{e}^f-1)}\right)>0\right\}
    \end{gathered}
    \end{equation}
The above domain can be observed on the right side of Figure \ref{fig:Sol_plane}.
\end{rem}

\subsection{Translation-like triangular surface in $\NIL$}
Based on the experience gained in the previous subsection, it is easy to develop the translation-like triangular surface in the remaining Thurston geometries. Therefore, we will not provide such a detailed explanation of the calculations. 

Let $P_1=(1,0,0,0)$, $P_2=(1,a,b,c)$ and $P_3=(1,d,e,f)$ be points in $\NIL$, and let $P=(1,x,y,z)$ be an arbitrary point on the translation-like triangular surface of $S^{\NIL,t}_{P_1,P_2,P_3}$. Pulling back $P$ to the origin, we obtain:
\begin{equation}
\begin{gathered}
\bT^{-1}_{P}(P_1)=P^1=(1,-x,-y,x\,y-z)\\
\bT^{-1}_{P}(P_2)=P^2=(1,a-x,b-y,c-z-b\,x+x\,y)\\
\bT^{-1}_{P}(P_3)=P^3=(1,d-x,e-y,f-z-e\,x+x\,y)\   
\end{gathered}
\label{nilbetoltak}
\end{equation}

The scaled tangent vector at the origin of the translation curve from the origin to any $P(1,x,y,z)$ can be calculated by the following formula, considering (\ref{2.14}):
\begin{equation}
   \bt(x,y,z)_\NIL:=d^t(E_0,P)\cdot\bt_e(x,y,z)=\left(x,y, z-\dfrac{x\,y}{2}\right)
    \label{nilerintov}
\end{equation}

Applying (\ref{nilerintov}) to every translated point (\ref{nilbetoltak}), we obtain $\bt(P_i)$ tangent vectors that are coplanar iff their triple product is 0. The result can be simplified to get an explicit result directly:
\begin{lem}
    Let $P_1=(1,0,0,0),$ $P_2=(1,a,b,c)$ and $P_3=(1,d,e,f)$  be points in the projective model of $\NIL$ geometry. Then the $S^{\NIL,t}_{P_1,P_2,P_3}$ translation-like triangular surface of $P_1,$ $P_2$ and $P_3$ has the equation:
\begin{equation}
    \begin{gathered}
       z=\dfrac{x\,y}{2}+
       \left(\dfrac{b\,d\,e-a\,b\,e+2c\,e-2b\,f}{2(a\,e-b\,d)}\right)\,x+
       \left(\dfrac{a\,b\,d-a\,d\,e+2a\,f-2c\,d}{2(a\,e-b\,d)}\right)\,y
    \end{gathered}
    \label{nilsik}
\end{equation}
\end{lem}

\begin{figure}[h]
    \centering
    \includegraphics[width=0.45\linewidth]{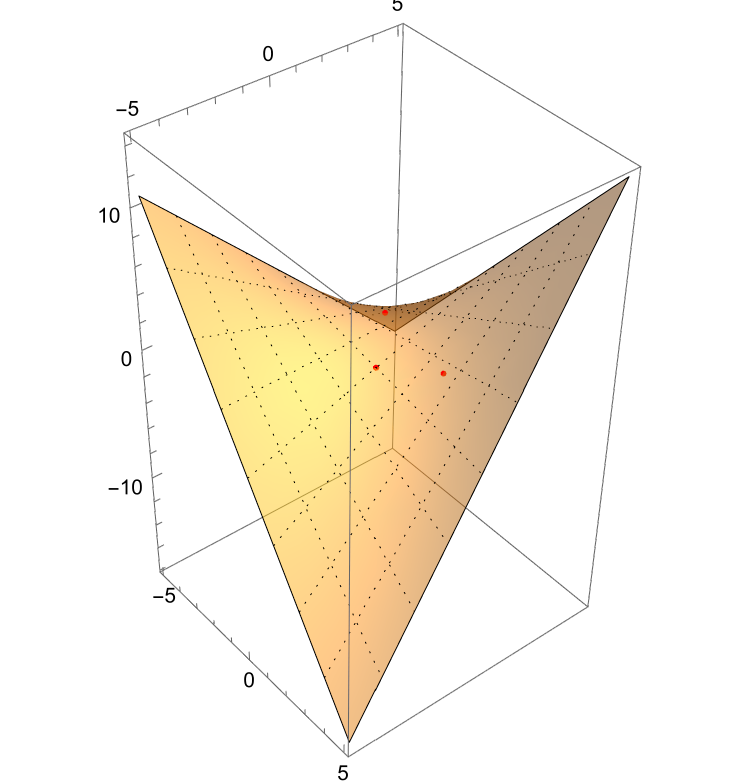}
    \includegraphics[width=0.45\linewidth]{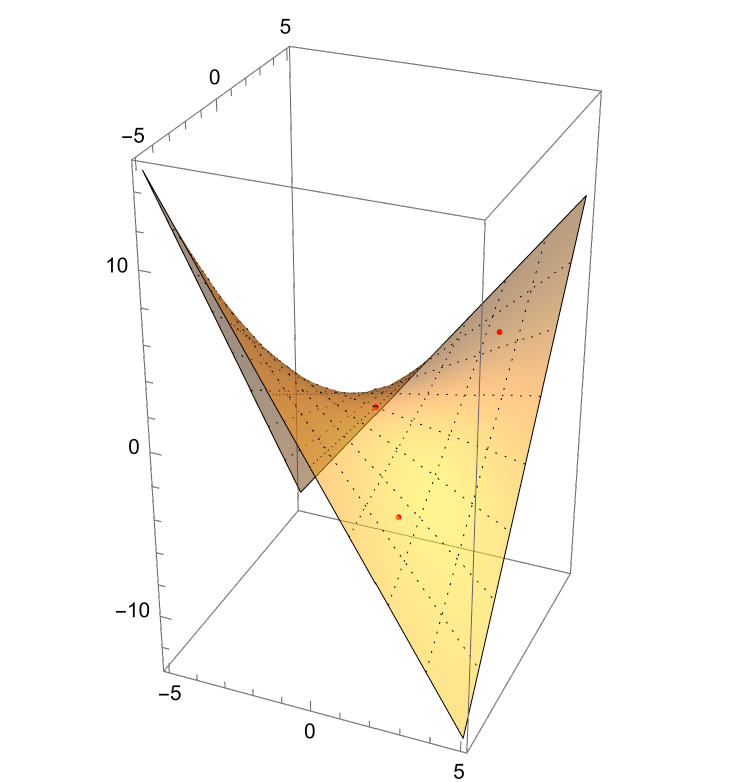}
    \caption{Translation-like triangular surface in $\NIL$ for points $P_1=(1,0,0,0),$ $P_2=(1,2,1,1), P_3=(1,-2,2,0)$ (left) and for points $P_1=(1,0,0,0),$ $P_2=(1,2,-3,-3), P_3=(1,3,3,3)$ (right).}
    \label{fig:nilplane}
\end{figure}

Some examples can be seen on Figure \ref{fig:nilplane}.
\begin{rem}
    \begin{enumerate}
        \item We excluded the $a\,e-b\,d=0$ case, while we transformed the result to an explicit expression. In that case, the result is a plane, through the $z$ axis.
        \item The graph of the above function is a hyperbolic paraboloid, \textit{i. e.} the result is a doubly ruled surface. So we can see that in the so called "twisted" geometries $(\SOL,\NIL,\SLR)$ the translation-like triangular surfaces are ruled surfaces.
    \end{enumerate}
\end{rem}
\begin{lem}
    Let $P_4$ be a point on $S^{\NIL,t}_{P_1,P_2,P_3}$ translation-like triangular surface. Then $S^{\NIL,t}_{P_1,P_2,P_3}=S^{\NIL,t}_{P_1,P_2,P_4}$ i. e. the definition of the translation-like triangular surface is transitive.
\end{lem}
\textsc{Proof}: The steps of the proof are the same as in the $\SOL$ case. $\square$

\subsection{Translation-like triangular surfaces in $\SXR$ and $\HXR$}

Now, we are interested in $\SXR$ and $\HXR$ spaces, at the same time. This is primarily due to the fact that a significant symmetry can be observed in both the calculations and the results. We will examine the two spaces simultaneously.

Let $P_1=(1,1,0,0)$, $P_2=(1,a,b,c)$ and $P_3=(1,d,e,f)$ be points in $\SXR $ or $\HXR$, and let $P=(1,x,y,z)$ be an arbitrary point on the translation-like triangular surface of $S^{\SXR,t}_{P_1,P_2,P_3}$ or $S^{\HXR,t}_{P_1,P_2,P_3}$. First, we need the matrix of the transformation, that pulls back $P$ to the origin of the model. A detailed explanation can be found in \cite{Cs23} or in \cite{Sz20} about this calculation, thus here we just provide the results. We remind the dear reader that $\pm$ is $+$ for $\SXR$ and $-$ for $\HXR:$
\footnotesize
\begin{equation}
\begin{gathered}
\bT_P^{-1}=
\begin{pmatrix}
1 & 0 & 0 & 0 \\
0 & \dfrac{x}{x^2\pm(y^2+z^2)} & -\dfrac{y}{x^2\pm(y^2+z^2)} & -\dfrac{z}{x^2\pm(y^2+z^2)} \\
0 & \pm\dfrac{y}{x^2\pm(y^2+z^2)} & \dfrac{xy^2+z^2\sqrt{x^2\pm(y^2+z^2)}}{(x^2\pm(y^2+z^2))(y^2+z^2)} & \dfrac{xyz-yz^2\sqrt{x^2\pm(y^2+z^2)}}{(x^2\pm(y^2+z^2))(y^2+z^2)} \\
0 & \pm\dfrac{z}{x^2\pm(y^2+z^2)} &  \dfrac{xyz-yz^2\sqrt{x^2\pm(y^2+z^2)}}{(x^2\pm(y^2+z^2))(y^2+z^2)} & \dfrac{xz^2+y^2\sqrt{x^2\pm(y^2+z^2)}}{(x^2\pm(y^2+z^2))(y^2+z^2)} 
\end{pmatrix}
\label{shxr_transz}
\end{gathered}
\end{equation}

\begin{equation}
\begin{gathered}
\bT_{P}^{-1}(P)=(1,1,0,0); \\
\bT_{P}^{-1}(P_1)=\left(1,\dfrac{x}{x^2\pm(y^2+z^2)},-\dfrac{y}{x^2\pm(y^2+z^2)},-\dfrac{z}{x^2\pm(y^2+z^2)}\right);\\
\bT_{P}^{-1}(P_2)=\begin{pmatrix} 1\\
\dfrac{ax\pm(by+cz)}{x^2\pm(y^2+z^2)}\\
\dfrac{zy(cx-az)-y^2(ay-bx)+z(bz-cy)\sqrt{x^2\pm(y^2+z^2)}}{(x^2\pm(y^2+z^2))(y^2+z^2)} \\ 
\dfrac{z^2(cx-az)-yz(ay-bx)-y(bz-cy)\sqrt{x^2\pm(y^2+z^2)}}{(x^2\pm(y^2+z^2))(y^2+z^2)} 
\end{pmatrix}^T;\\
\bT_{P}^{-1}(P_3)=\begin{pmatrix} 1\\
\dfrac{dx\pm(ey+fz)}{x^2\pm(y^2+z^2)}\\
\dfrac{zy(fx-dz)-y^2(dy-ex)+z(ez-fy)\sqrt{x^2\pm(y^2+z^2)}}{(x^2\pm(y^2+z^2))(y^2+z^2)} \\ 
\dfrac{z^2(fx-dz)-yz(dy-ex)-y(ez-fy)\sqrt{x^2\pm(y^2+z^2)}}{(x^2\pm(y^2+z^2))(y^2+z^2)} 
\end{pmatrix}^T
\label{shxreltoltak}
\end{gathered}
\end{equation}
\normalsize

We need the tangent vectors of the translation curves running into the points $\bT^{-1}_{P}(P_i)$ $(i=1,2,3)$ from $\bT^{-1}_{P}(P)=E_0.$ It is not necessary to determine the exact value of the parameters $u,$ $v,$ $\tau$, it is enough to evaluate the vector $\bt$ (see \ref{2.8}). 

\begin{lem} 
Let $(1,x,y,z)$ $(x,y,z\in\mathbb{R}$ and $x^2+y^2+z^2>0$ in $\SXR$; $x^2-y^2-z^2>0,\ x>0$ in $\HXR)$ be the homogeneous coordinates of a point $P\in\SXR$ or $P\in\HXR.$ Then the translation curve, drawn to $P$ from $E_0=(1,1,0,0)$ has the following tangent in $E_0$:
\footnotesize
\begin{equation}
\tau\cdot\bt_P=\left(\dfrac{1}{2}\ln(x^2\pm(y^2+z^2)),
\dfrac{y\,\mathrm{arcC}\left(\dfrac{x}{\sqrt{x^2\pm(y^2+z^2)}}\right)}{\sqrt{y^2+z^2}},
\dfrac{z\,\mathrm{arcC}\left(\dfrac{x}{\sqrt{x^2\pm(y^2+z^2)}}\right)}{\sqrt{y^2+z^2}}\right)
\end{equation}
\normalsize
where $\pm$ and $\mathrm{arcC}(x)$ is $+$ and $\arccos(x)$ for $\SXR,$ $-$ and $\mathrm{arccosh}(x)$ for $\HXR,$ and $\tau$ is the distance of $P$ and $E_0.$ $\square$ 
\label{tan_lem}
\end{lem}

As above, after determining the tangents at the origin, we are able to consider the triple product of the tangents to obtain the following lemma.  Due to the length of the formula, the result is only presented in its seriously simplified form, introducing some notations (see Theorem \ref{shxrsik}). During the calculations we used dot and cross products and the spherical/hyperbolic law of sine and cosine for sides.

\begin{thm}
Let $P_1=(1,1,0,0),$ $P_2=(1,a,b,c)$ and $ P_3=(1,d,e,f)$ be given points in $\SXR$ or $\HXR$ and $P=(1,x,y,z)$ be an arbitrary point. If $\bp_1=(1,0,0);$ $\bp_2=(a,b,c);$ $\bp_3=(d,e,f);$ $\bp=(x,y,z);$ the projected image of $P_1,$ $P_2,$ $P_3$ and $P$ onto the unit surface (sphere in $\SXR$ and hyperboloid in $\HXR$) of the geometry along the corresponding fibre lines are $P_1'=P_1$, $P_2'$, $P_3'$  
and $P';$ 
$d_1,$ $d_2$ and $d_3$ are the distances of $P'$ to $P_1',$ to $P_2'$ and to $P_3'$ and $\gamma_1=P_2'P'P_3'\angle,$ $\gamma_2=P_3'P'P_1'\angle,$ $\gamma_3=P_1'P'P_2'\angle$ are directed angles; then $P$ lies on the $S^{\SXR,t}_{P_1,P_2,P_3}$ or $S^{\HXR,t}_{P_1,P_2,P_3}$ translation-like triangular surface of $P_1,$ $P_2$ and $P_3$ if
\begin{equation}
d_1d_2\sin(\gamma_3)\ln\frac{|\bp_3|}{|\bp|}+
d_2d_3\sin(\gamma_1)\ln\frac{|\bp_1|}{|\bp|}+
d_3d_1\sin(\gamma_2)\ln\frac{|\bp_2|}{|\bp|}=0,
\label{eq:shxrsik}
\end{equation} 
where $|\cdot|$ is the Euclidean norm. $\ \ \square$
\label{shxrsik}
\end{thm}

\begin{rem}
    Fixing $P'$ point on the surface of the unit surface of the geometry will determine $d_i$ and $\gamma_i$ $(i=1,2,3)$. Therefore $\ln|\bp|$ can be expressed from $(\ref{eq:shxrsik}):$
    \begin{equation}
        \ln|\bp|=\dfrac{d_1d_2\sin(\gamma_3)\ln|\bp_3|+
d_2d_3\sin(\gamma_1)\ln|\bp_1|+
d_3d_1\sin(\gamma_2)\ln|\bp_2|}{d_1d_2\sin(\gamma_3)+
d_2d_3\sin(\gamma_1)+
d_3d_1\sin(\gamma_2)},
    \end{equation}
which means that along every fibre line, the translation-like triangular surface has exactly 1 point.    
\end{rem}
\begin{figure}
    \centering
    \includegraphics[width=0.45\linewidth]{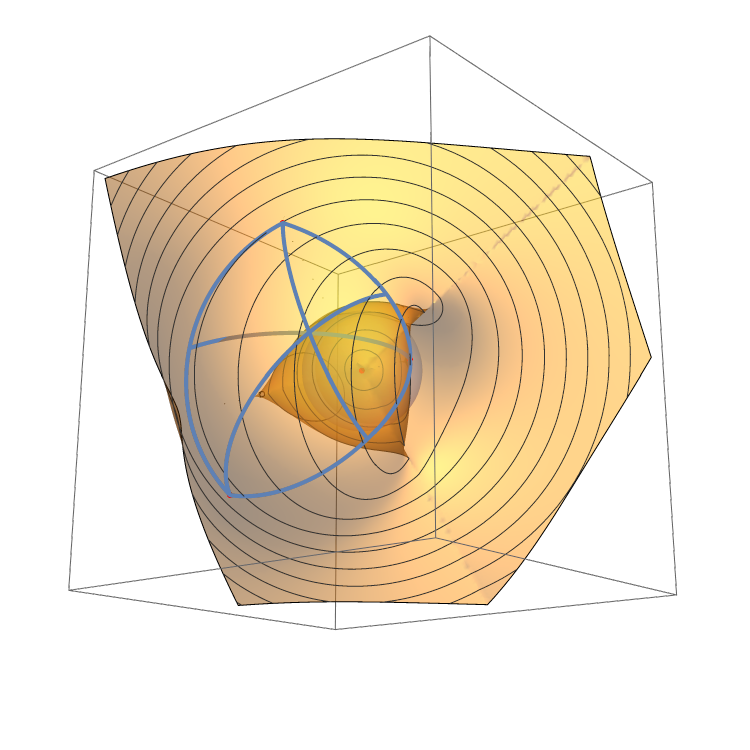}    \includegraphics[width=0.45\linewidth]{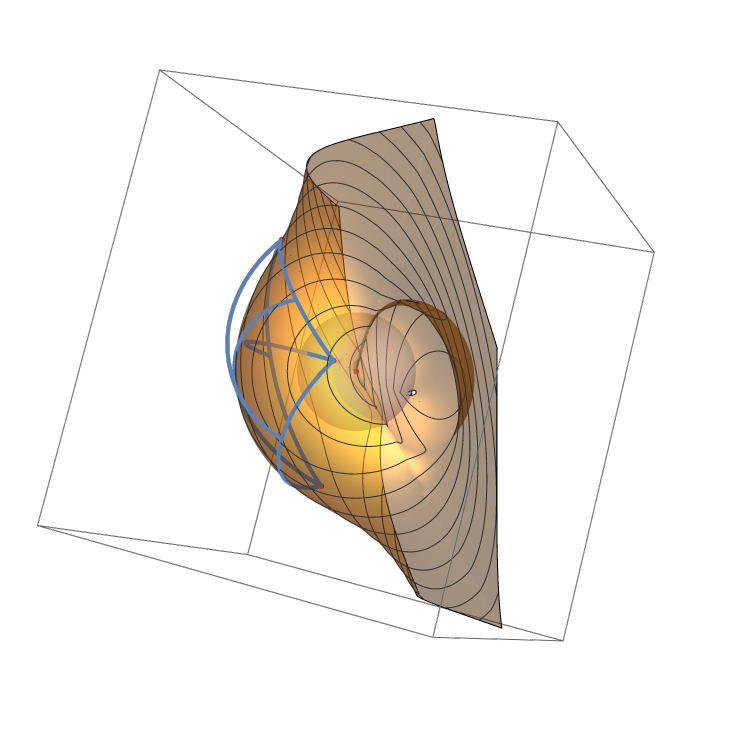}
    \caption{Translation-like triangular surface in $\SXR$ for $P_1=(1,1,0,0)$, $P_2=(1,0,-2,2)$ and $P_3=(1,-2,-1,-2)$}
    \label{fig:s2xr_plane}
\end{figure}

\begin{figure}
    \centering
    \includegraphics[width=0.45\linewidth]{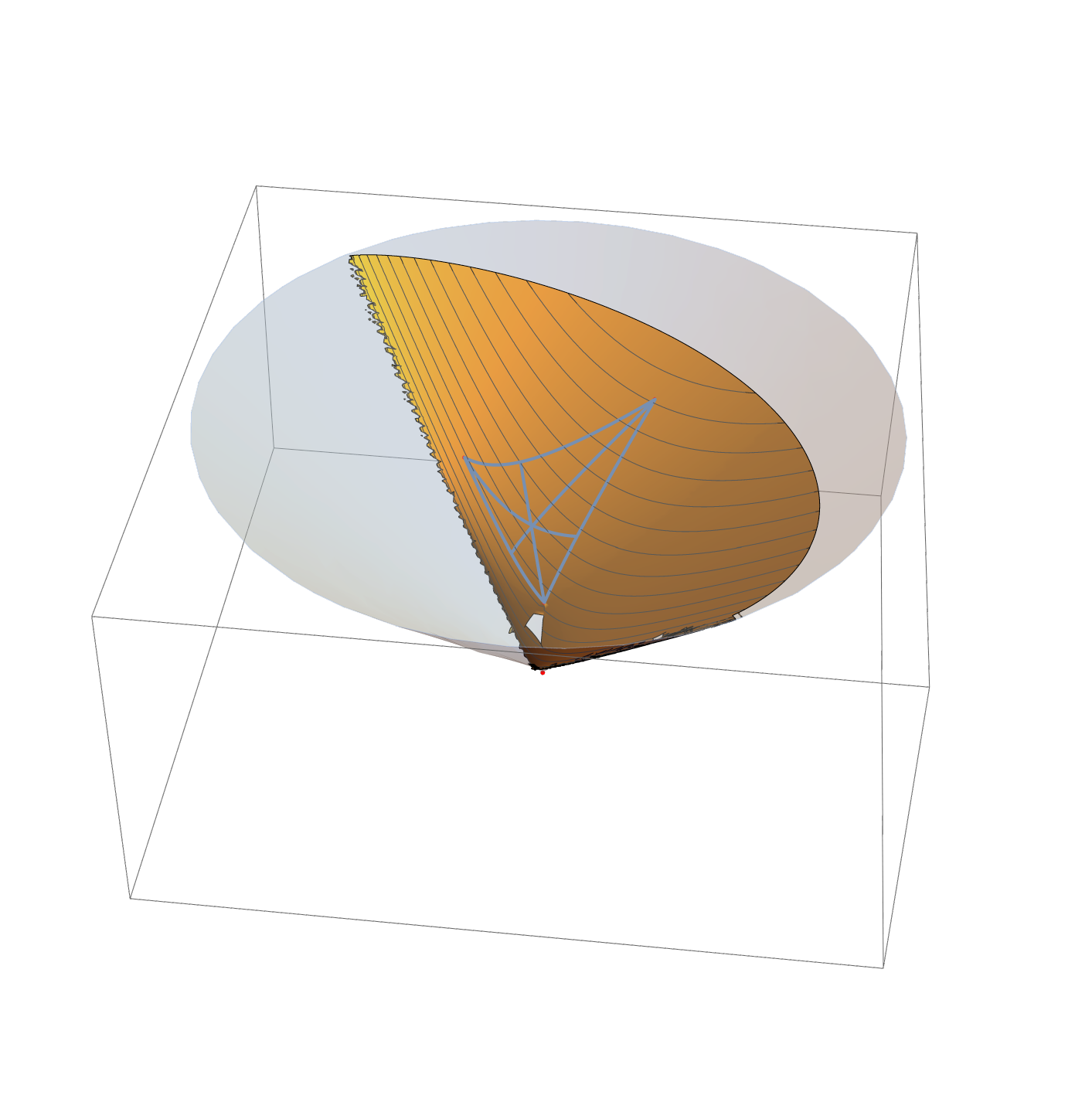}    \includegraphics[width=0.45\linewidth]{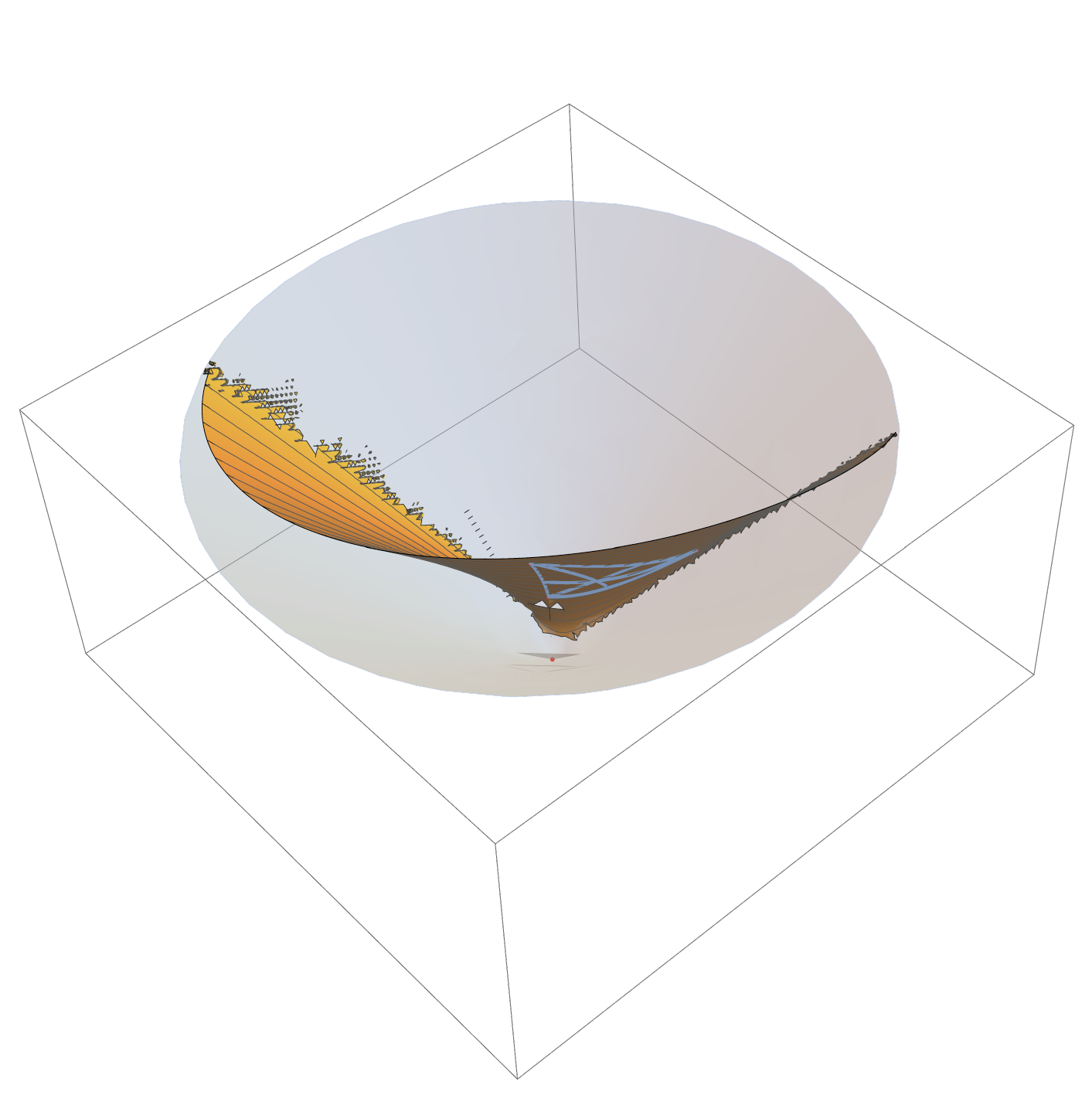}
    \caption{Translation-like triangular surface in $\HXR$ for $P_1=(1,1,0,0)$, $P_2=(1,3,1,1)$ and $P_3=(1,1.9,-1,1.2)$}
    \label{fig:h2xr_plane}
\end{figure}

The yellow surface on Figure \ref{fig:s2xr_plane}-\ref{fig:h2xr_plane} shows an example for translation-like triangular surface in $\SXR$ and $\HXR$ geometries. It is interesting, but not new (see \cite{Sz21}) that although the (blue) translation curves between the given (red) points are on the translation-like triangular surface, e. g. the medians of the triangle are neither lie on the surface nor intersect each other (see the right side of Figure \ref{fig:s2xr_plane}). This is the key of the proof in the following lemma:
\begin{lem}
    Let $P_4$ be a point on $S^{X,t}_{P_1,P_2,P_3}$ translation-like triangular surface, where $X$ is either $\SXR$ or $\HXR$. Then $S^{X,t}_{P_1,P_2,P_3}\neq S^{X,t}_{P_1,P_2,P_4}$ i. e. the definition of the translation-like triangular surface is \textbf{not} transitive.
\end{lem}
\textsc{Proof}: Let $P_4$ be the midpoint of $P_1$ and $P_3$, then $S^{X,t}_{P_1,P_2,P_4}$ contains the translation curve from $P_4$ to $P_2$, but $S^{X,t}_{P_1,P_2,P_3}$ doesn't. $\square$
%

\end{document}